\input ./style/arxiv-general.cfg
\documentclass[aop,MSNbibl,seceqn,dvips]{arximspdf}
\makeatletter
   \@ifpackageloaded{graphicx}{}{\usepackage{graphicx}}
\makeatother

\doi{10.1214/14-AOP979} 
\volume{44}
\issue{1}
\pubyear{2016}
\firstpage{479}
\lastpage{520}
\docsubty{FLA}

\makeatletter
\newcommand{\rrvert}{\vert}
\newcommand{\llvert}{\vert}
\newcommand{\eqref}[1]{(\ref{#1})}
\newtheorem{thmm}{Theorem}[section]
\newtheorem{lemma}[thmm]{Lemma}
\newtheorem{prop}[thmm]{Proposition}
\newtheorem{cor}[thmm]{Corollary}
\newproclaim{defn}[thmm]{Definition}
\newproclaim{examp}[thmm]{Example}
\newproclaim{rmk}[thmm]{Remark}

\newcommand{\R}{\mathbb{R}}
\newcommand{\Z}{\mathbb{Z}}
\newcommand{\N}{\mathbb{N}}
\newcommand{\E}{\mathbb{E}}
\newcommand{\BBW}{\mathbb{W}}
\newcommand{\BBM}{\mathbb{M}}
\newcommand{\BBU}{\mathbb{U}}
\newcommand{\Cov}{\operatorname{Cov}}
\newcommand{\essinf}{\operatorname{ess\,inf}}
\newcommand{\tildev}{\tilde{v}}
\newcommand{\tildeW}{\widetilde{W}}
\newcommand{\tildeBBW}{\widetilde{\mathbb W}}
\makeatother

\begin{document}
\begin{frontmatter}

\title{Smooth approximation of stochastic differential~equations}
\runtitle{Smooth approximation of SDEs}

\begin{aug}
\author[A]{\fnms{David}~\snm{Kelly}\ead[label=e1]{dtkelly@cims.nyu.edu}\thanksref{T1}}
\and
\author[B]{\fnms{Ian}~\snm{Melbourne}\corref{}\ead[label=e2]{i.melbourne@warwick.ac.uk}\thanksref{T2}}
\affiliation{University of North Carolina and University of Warwick}
\address[A]{Courant Institute of Mathematical Sciences\\
New York University\\
New York, New York 10012-1185\\
USA\\
\printead{e1}}
\address[B]{Mathematics Institute\\
University of Warwick\\
Coventry CV4 7AL\\
United Kingdom\\
\printead{e2}}
\thankstext{T1}{Supported in part by ONR Grant N00014-12-1-0257.}
\thankstext{T2}{Supported in part by European Advanced Grant
StochExtHomog (ERC AdG 320977).}

\runauthor{D. Kelly and I. Melbourne}
\end{aug}

\received{\smonth{3} \syear{2014}}
\revised{\smonth{10} \syear{2014}}

\begin{abstract}
Consider an It\^o process $X$ satisfying the stochastic differential equation
$dX=a(X) \,dt+b(X) \,dW$ where $a,b$ are smooth and $W$ is a
multidimensional Brownian motion.
Suppose that $W_n$ has smooth sample paths and that $W_n$ converges
weakly to $W$. A central question in stochastic analysis is to
understand the limiting behavior of solutions $X_n$ to the ordinary
differential equation $dX_n=a(X_n) \,dt+b(X_n) \,dW_n$.

The classical Wong--Zakai theorem gives sufficient conditions under
which $X_n$ converges weakly to $X$ provided that the stochastic
integral $\int b(X) \,dW$ is given the Stratonovich interpretation. The
sufficient conditions are automatic in one dimension, but in higher
dimensions the correct interpretation of
$\int b(X) \,dW$ depends sensitively on how the smooth approximation
$W_n$ is chosen.

In applications, a natural class of smooth approximations arise by setting
$W_n(t)=n^{-1/2}\int_0^{nt}v\circ\phi_s \,ds$ where $\phi_t$ is a flow
(generated, e.g., by an ordinary differential equation) and $v$
is a mean zero observable.
Under mild conditions on $\phi_t$, we give a definitive answer to the
interpretation question for the stochastic integral $\int b(X) \,dW$.
Our theory applies to Anosov or Axiom~A flows $\phi_t$, as well as to a
large class of nonuniformly hyperbolic flows (including the one defined
by the well-known Lorenz equations) and our main results do not require
any mixing assumptions on $\phi_t$.

The methods used in this paper are a combination of rough path theory
and smooth ergodic theory.
\end{abstract}

%

%
\end{frontmatter}

\section{Introduction}
\label{sec-intro}

Let $X$ be a $d$-dimensional It\^o process defined by a stochastic
differential equation (SDE) of the form
%
\begin{equation}
\label{eq-IP} dX=a(X) \,dt+b(X) \,dW,
\end{equation}
where $a\dvtx\R^d\to\R^d$ is $C^{1+}$, $b\dvtx\R^d\to\R
^{d\times e}$ is $C^{2+}$,
and $W$ is an $e$-dimensional Brownian motion with $e\times e$-dimensional
covariance matrix $\Sigma$.

Given a sequence of $e$-dimensional processes $W_n$ with smooth sample
paths, we consider the sequence of ordinary differential equations (ODEs)
%
\begin{equation}
\label{eq-ODE} dX_n=a(X_n) \,dt+b(X_n)
\,dW_n,
\end{equation}
where $dW_n=\dot W_n \,dt$.
We suppose that an initial condition $\xi\in\R^d$ is fixed throughout
and consider solutions $X$ and $X_n$ satisfying $X(0)=X_n(0)=\xi$.

Let $T>0$. The sequence $W_n$ is said to satisfy the \emph{weak
invariance principle (WIP)}
if $W_n\to_w W$ in $C([0,T],\R^e)$. Assuming the WIP, a central
question in stochastic analysis is to determine whether $X_n\to_w X$ in
$C([0,T],\R^d)$ for a suitable interpretation of the stochastic
integral $\int b(X) \,dW$ implicit in \eqref{eq-IP}.
The Wong--Zakai theorem \cite{WongZakai65} gives general conditions
under which convergence holds with the Stratonovich interpretation for
the stochastic integral.
These conditions are automatically satisfied in the one-dimensional
case $d=e=1$, but may fail in higher dimensions. See also Sussmann
\cite
{Sussmann78}.
In two dimensions, McShane~\cite{McShane72} gave the first counterexamples,
and Sussmann \cite{Sussmann91} provided numerous further counterexamples.

From now on, we replace \eqref{eq-IP} by the SDE
%
\begin{equation}
\label{eq-SDE} dX=a(X) \,dt+b(X)\ast dW,
\end{equation}
to emphasize the issue with the interpretation of the stochastic integral.
General principles suggest that the limiting stochastic integral should
be Stratonovich modified by an antisymmetric drift term:
\[
b(X)\ast \,dW=b(X)\circ \,dW+\frac{1}2\sum_{\alpha,\beta,\gamma}D^{\beta
\gamma}
\partial^\alpha b^{\beta}(X) b^{ \alpha\gamma}(X) \,dt.
\]
Here, and throughout the paper, we sum over $1\le\alpha\le d$, $1\le
\beta,\gamma\le e$, and $b^{\alpha\gamma}$ and $b^\beta$ denote
the ($\alpha,\gamma$)th entry and
$\beta$th column, respectively, of $b$. Moreover, $\{D^{\beta\gamma
}\}
$ is an antisymmetric
matrix that is to be determined.
[Hence, an alternative to \eqref{eq-SDE} would be to consider
$dX=\tilde a(X) \,dt+b(X)\circ dW$ with the emphasis on determining the
correct drift term $\tilde a$.]

In applications, smooth processes $W_n$ that approximate Brownian
motion arise naturally from differential equations as
follows \cite
{GivonKupfermanStuart04,Huisinga03,MS11,Papanicolaou74,PavliotisStuart}.
Let $\phi_t\dvtx M\to M$ be a smooth flow on a finite-dimensional manifold
$M$ preserving an ergodic measure $\nu$ and
let $v\dvtx M\to\R^e$ be a smooth observable with $\int_M v \,d\nu=0$.
Define
%
\begin{equation}
\label{eq-W} W_n(t)=n^{-1/2}\int_0^{nt}
v\circ\phi_s \,ds.
\end{equation}
For large classes of uniformly and nonuniformly hyperbolic flows \mbox{\cite
{DenkerPhilipp84,MN05,MN09,Gouezel10}}, it can be shown that $W_n$
satisfies the WIP.
In this paper, we consider such flows, and give a definitive answer to
the question of how to correctly interpret the stochastic integral
$\int b(X)\ast \,dW$ in order to ensure that $X_n\to_w X$.

\textit{An important special case}.
Let $d=e=2$ and take $a\equiv0$, $b(x_1,x_2)= \bigl(
{1 \enskip 0 \atop 0 \enskip x_1} \bigr)$.
The ODE \eqref{eq-ODE} becomes
\[
dX_n^1=dW_n^1,\qquad
dX_n^2 = X_n^1
\,dW_n^2,
\]
so with the initial condition\vspace*{1pt} $\xi=0$ we obtain $X_n^1\equiv W_n^1$ and
$X_n^2(t)=\int_0^t W_n^1 \,dW_n^2$.
Weak convergence of $W_n$ to $W$ does not determine the weak limit of
$\int_0^t W_n^1 \,dW_n^2$.
However,
according to rough path theory \cite{Lyons98}, this is the key
obstruction to solving the central problem in this paper.
Generally, define the family of smooth processes $\BBW_n\in
C([0,\infty
),\R^{e\times e})$,
%
\begin{equation}
\label{eq-WW} \BBW_n^{\beta\gamma}(t)=\int_0^t
W_n^\beta \,dW_n^\gamma,\qquad 1\le \beta,
\gamma\le e.
\end{equation}
The theory of rough paths implies that under some mild moment estimates,
the weak limit of $(W_n,\BBW_n)$ determines the weak limit of $X_n$
in \eqref{eq-ODE} and the correct interpretation for the stochastic
integral in \eqref{eq-SDE}.

Hence, a large part of this paper is dedicated to proving an \emph
{iterated WIP} for the pair $(W_n,\BBW_n)$.

\textit{Anosov and Axiom A flows}.
One well-known class of flows to which our results apply is given by
the Axiom A (uniformly hyperbolic) flows introduced by Smale \cite{Smale67}.
This includes Anosov flows \cite{Anosov67}. We do not give the precise
definitions, since they are not needed for understanding the paper, but
a rough description is as follows. (See \cite{Ruelle78,Bowen75,Sinai72}
for more details.)

Let $\phi_t\dvtx M\to M$ be a $C^2$ flow defined on a compact manifold $M$.
A flow-invariant subset $\Omega\subset M$ is \emph{uniformly hyperbolic}
if for all
$x\in\Omega$ there exists a $D\phi_t$-invariant splitting transverse to
the flow into uniformly contracting and expanding directions.
The flow is \emph{Anosov} if the whole of $M$ is uniformly hyperbolic.
More generally, an \emph{Axiom A} flow is characterised by the property
that the dynamics decomposes into finitely many hyperbolic equilibria
and finitely many uniformly hyperbolic subsets $\Omega_1,\ldots
,\Omega
_k$, called \emph{hyperbolic basic sets}, such that the flow on each
$\Omega_i$ is transitive (there is a dense orbit).

If $\Omega$ is a hyperbolic basic set, there is a unique $\phi
_t$-invariant ergodic probability measure (called an \emph{equilibrium
measure}) associated to
each H\"older function on $\Omega$. [In the special case that $\Omega$
is an attractor, there is a distinguished equilibrium measure called
the physical measure or SRB measure (after Sinai, Ruelle, Bowen).]

In the remainder of the \hyperref[sec-intro]{Introduction}, we assume
that $\Omega$ is a
hyperbolic basic set
with equilibrium measure $\nu$ (corresponding to a H\"older potential).
We exclude the trivial case where $\Omega$ consists of a single
periodic orbit.

We can now state our main results.
For $u\dvtx\Omega\to\R^q$, we define $\E_\nu(u)\in\R^q$
and $\Cov_\nu(u)\in\R^{q\times q}$ by setting
$\E_\nu(u)=\int_\Omega u \,d\nu$ and
$\Cov_\nu^{\beta\gamma}(u)=\E_\nu(u^\beta u^\gamma)-\E_\nu
(u^\beta)\E
_\nu(u^\gamma)$.

\begin{thmm}[(Iterated WIP)] \label{thmm-WW}
Suppose that $\Omega\subset M$ is a hyperbolic basic set with
equilibrium measure $\nu$ and that $v\dvtx\Omega\to\R^e$ is H\"
older with
$\int_\Omega v \,d\nu=0$.
Define $W_n$ and $\BBW_n$ as in \eqref{eq-W} and \eqref{eq-WW}. Then:
\begin{enumerate}[(a)]
\item[(a)]
$(W_n,\BBW_n)\to_w (W,\BBW)$ in
$C([0,\infty),\R^e\times\R^{e\times e})$ as $n\to\infty$, where:
\begin{enumerate}[(ii)]
\item[(i)]
$W$ is an $e$-dimensional
Brownian motion with covariance matrix $\Sigma=\Cov(W(1))=\lim_{n\to
\infty}\Cov_\nu(W_n(1))$.
\item[(ii)]
$\BBW^{\beta\gamma}(t)=\int_0^t W^\beta\circ dW^\gamma+ \frac{1}2
D^{\beta\gamma}t$
where $D=2\lim_{n\to\infty}\E_\nu(\BBW_n(1))-\Sigma$.
\end{enumerate}
\item[(b)] If in addition the integral
$\int_0^\infty\int_{\Omega}v^\beta v^\gamma\circ\phi_t \,dt$ exists
for all $\beta,\gamma$, then
\[
\Sigma^{\beta\gamma}= \int_0^\infty\int
_\Omega\bigl(v^\beta v^\gamma\circ
\phi_t+ v^\gamma v^\beta\circ\phi_t
\bigr) \,d\nu \,dt
\]
and
\[
D^{\beta\gamma}= \int_0^\infty\int
_\Omega\bigl(v^\beta v^\gamma\circ
\phi_t- v^\gamma v^\beta\circ\phi_t
\bigr) \,d\nu \,dt.
\]
\end{enumerate}
\end{thmm}

\begin{thmm}[(Convergence to SDE)] \label{thmm-approx}
Suppose that $\Omega\subset M$ is a hyperbolic basic set with
equilibrium measure $\nu$ and that $v\dvtx X\to\R^e$ is H\"older
with $\int_\Omega v \,d\nu=0$.
Let $W_n$, $W$ and $D$ be as in Theorem~\ref{thmm-WW}.
Let $a\dvtx\R^d\to\R^d$ be $C^{1+}$ and $b\dvtx\R^d\to\R
^{d\times e}$ be
$C^{2+}$, and
define $X_n$ to be the solution of the ODE \eqref{eq-ODE} with
$X_n(0)=\xi$.

Then
$X_n\to_w X$ in $C([0,\infty),\R^d)$ as $n\to\infty$, where $X$
satisfies the SDE
\[
dX= \biggl\{a(X)+\frac{1}2\sum_{\alpha,\beta,\gamma}D^{\beta\gamma
}
\partial ^\alpha b^{\beta}(X) b^{ \alpha\gamma}(X) \biggr\} \,dt+b(X)
\circ dW,\qquad X(0)=\xi.
\]
\end{thmm}

\textit{Mixing assumptions on the flow}.
The only place where we use mixing assumptions on the flow is in
Theorem~\ref{thmm-WW}(b) to obtain closed form expressions for the
diffusion and drift coefficients $\Sigma$ and $D$. In general, these
integrals need
not converge for Axiom A flows even when $v$ is $C^\infty$.

Dolgopyat \cite{Dolgopyat98a} proved exponential decay of correlations
for H\"older observables $v$ of certain Anosov flows, including
geodesic flows on compact negatively curved surfaces. This was extended
by Liverani \cite{Liverani04} to Anosov flows with a contact structure,
including the case of geodesic flows in all dimensions. Theorem~\ref
{thmm-WW}(b) holds for the flows considered
in \cite{Dolgopyat98a,Liverani04}. Nevertheless, for typical Anosov
flows, the extra condition in
Theorem~\ref{thmm-WW}(b) is not known to hold for H\"older observables.

Dolgopyat \cite{Dolgopyat98b} introduced the weaker notion of \emph
{rapid mixing}, namely decay of correlations at an arbitrary polynomial
rate, and proved that typical Axiom A flows enjoy this property.
By \cite{FMT07}, an open and dense set of Axiom A flows are rapid
mixing. However, this theory applies only to observables $v$ that are
sufficiently smooth, and the degree of smoothness
is not readily computable. On the positive side, Theorem~\ref
{thmm-WW}(b) holds for typical Axiom A flows provided $v$ is $C^\infty$.

In the absence of a good theory of mixing for flows, we have chosen (as
in~\cite{MS11}) to develop our theory in such a way that the dependence
on mixing is minimized. Instead we rely on statistical properties of
flows, which is a relatively well-understood topic.

A more complicated closed form expression for $\Sigma$ and $D$ that
does not require mixing conditions on the flow can be found in
Corollary~\ref{cor-WWflow}.

\textit{Beyond uniform hyperbolicity}.
In this \hyperref[sec-intro]{Introduction}, for ease of exposition we
have chosen to focus on
the case of
uniformly hyperbolic flows (Anosov or Axiom~A).
However, our results hold for large classes of nonuniformly hyperbolic flows.
In particular, Young \cite{Young98} introduces a class of nonuniformly
hyperbolic diffeomorphisms, that includes uniformly hyperbolic
(Axiom A) diffeomorphisms, as well as H\'enon-like attractors \cite
{BenedicksYoung00}.
For flows with a Poincar\'e map that is nonuniformly hyperbolic in the
sense of \cite{Young98}, Theorems \ref{thmm-WW} and \ref
{thmm-approx} go
through unchanged.

The nonuniformly hyperbolic diffeomorphisms in \cite{Young98} (but not
necessarily the corresponding flows) have exponential decay of
correlations for H\"older observables. Young \cite{Young99} considers
nonuniformly hyperbolic diffeomorphisms with subexponential decay of
correlations. Many of our results go through for flows with a Poincar\'
e map that is nonuniformly hyperbolic in the more general sense of~\cite
{Young99}.
In particular, our results are valid for the classical Lorenz equations.

These extensions are discussed at length in Section~\ref{sec-gen}.

\textit{Structure of the proofs.}
In the smooth ergodic theory literature, there are numerous results on
the WIP where $W_n\to_w W$. Usually such results are obtained first for
processes $W_n$ arising from a discrete time dynamical system.
Results for flows are then obtained as a corollary of the discrete time case,
see for example \cite{Ratner73,MT04,MN05,MN09,BunimSinaiChernov91,MZapp}.
Hence, it is natural to solve the discrete time analogue of Theorem~\ref
{thmm-WW} first before extending to continuous time.
This is the approach followed in this paper.
We first prove the discrete time iterated WIP, Theorem~\ref{thmm-WW2}
below. Then we derive the continuous time WIP, Theorem~\ref{thmm-WW}, as
a consequence, before obtaining Theorem~\ref{thmm-approx} using rough
path theory.
For completeness, we also state and prove the discrete time analogue of
Theorem~\ref{thmm-approx} (see Theorem~\ref{thmm-approx2} below), even
though this is not required for the proof of Theorem~\ref{thmm-approx}.

For the proof of the discrete time iterated WIP, it is convenient to
use the standard method of passing from invertible maps to
noninvertible maps. So we prove the iterated WIP first for
noninvertible maps, then for invertible maps, and finally for
continuous time systems.

\textit{Structure of the paper}.
The remainder of this paper is organized as follows.
Sections \ref{sec-discrete} to \ref{sec-NUH}
deal with the discrete time iterated WIP.
Section~\ref{sec-discrete} states our main results for discrete time.
In Section~\ref{sec-cohom}, we present a result on cohomological
invariance of weak limits of iterated processes. This result seems of
independent theoretical interest but in this paper it is used to
significantly simplify calculations.
In Sections \ref{sec-NUE} and \ref{sec-NUH}, we prove the iterated WIP
for discrete time systems that are noninvertible and invertible, respectively.

In Section~\ref{sec-flow}, we return to the case of continuous time and
prove a purely probabilistic result about lifting the iterated WIP from
discrete time to continuous time.
In Section~\ref{sec-moment}, we state and prove some moment estimates
that are required to apply rough path theory.
In Section~\ref{sec-flowpf}, we
prove the iterated WIP stated in Theorem~\ref{thmm-WW}.
Then in Section~\ref{sec-approx}, we prove Theorem~\ref{thmm-approx} and
its discrete time analogue.

In Section~\ref{sec-gen}, we discuss various generalizations of
our main results that go beyond the Axiom A case. In particular, we
consider large classes of systems that are nonuniformly hyperbolic in
the sense of \cite{Young98,Young99}.

We conclude this \hyperref[sec-intro]{Introduction} by mentioning
related work of
Dolgopyat \cite{Dolgopyat04}, Theorem~5 and \cite{Dolgopyat05},
Theorem~3(b).
These results, which rely on very different techniques from those
developed here, prove the analogue of Theorem~\ref{thmm-approx} for a class
of partially hyperbolic discrete time dynamical systems. The
intersection with our work consists of Anosov diffeomorphisms and
time-one maps of Anosov flows with better than summable decay of
correlations. As discussed above, our main results do not rely on
mixing for flows; only the formulas require mixing. Also, we consider
the entire Axiom A setting (including Smale horseshoes and flows
that possess a horseshoe in the Poincar\'e map) and our results apply
to systems that are nonuniformly hyperbolic in the sense of Young
(including H\'enon and Lorenz attractors).

\textit{Notation.}
As usual, we let $\int b(X) \,dW$ and $\int b(X)\circ dW$
denote the It\^o and Stratonovich integrals, respectively.

We use the ``big $O$'' and $\ll$ notation interchangeably, writing
$a_n=O(b_n)$ or $a_n\ll b_n$ if there is a constant $C>0$ such that
$a_n\le Cb_n$ for all $n\ge1$.

\section{Statement of the main results for discrete time}
\label{sec-discrete}

In this section, we state the discrete time analogues of our main
Theorems \ref{thmm-WW} and \ref{thmm-approx}.

Let $f\dvtx M\to M$ be a $C^2$ diffeomorphism defined on a compact
manifold $M$.
Again we focus on the case where $\Lambda\subset M$ is a (nontrivial)
hyperbolic basic set with equilibrium measure $\mu$. The definitions
are identical to those for Axiom A flows, with the simplification that
the direction tangent to the flow is absent.
(Hyperbolic basic sets are denoted throughout by $\Omega$ in the flow
case described in Section~\ref{sec-intro}
and by $\Lambda$ in the current discrete time setting. The analysis of
the flow case includes passing from the hyperbolic basic set $\Omega$
for the flow to a hyperbolic basic set $\Lambda$ for a suitable
Poincar\'e map; hence the need for distinct notation.)

We assume in this section that $\Lambda$ is mixing: $\lim_{n\to
\infty
}\int_\Lambda w_1 w_2\circ f^n \,d\mu=\int_\Lambda w_1 \,d\mu\int_\Lambda w_2 \,d\mu$ for all $w_1,w_2\in L^2$ (this assumption is
relaxed in
Section~\ref{sec-gen}).

Let $v\dvtx\Lambda\to\R^e$ be H\"older with $\int_\Lambda v \,d\mu=0$.
Define the cadlag processes $W_n\in D([0,\infty),\R^e)$,
$\BBW_n\in D([0,\infty),\R^{e\times e})$,
%
\begin{eqnarray}
\label{eq-WWW} W_n(t) & =&n^{-1/2}\sum
_{j=0}^{[nt]-1}v\circ f^j,
\nonumber
\\[-8pt]
\\[-8pt]
\nonumber
\BBW^{\beta\gamma}_n(t) &=&\int_0^t
W_n^\beta \,dW_n^\gamma=n^{-1}
\sum_{0\le i<j\le[nt]-1} v^\beta\circ f^i
v^\gamma\circ f^j.
\end{eqnarray}
Since our limiting processes have continuous sample paths, throughout
we use the sup-norm topology on $D([0,\infty),\R^e)$ unless otherwise stated.

\begin{thmm}[(Iterated WIP, discrete time)] \label{thmm-WW2}
Suppose that $\Lambda\subset M$ is a mixing hyperbolic basic set with
equilibrium measure $\mu$, and that $v\dvtx\Lambda\to\R^e$ is H\"
older with
$\int_\Lambda v \,d\mu=0$.
Define $W_n$ and $\BBW_n$ as in \eqref{eq-WWW}.
Then
$(W_n,\BBW_n)\to_w (W,\BBW)$ in $D([0,\infty),\R^e\times\R
^{e\times e})$
as $n\to\infty$, where:
\begin{longlist}[(ii)]
\item[(i)]
$W$ is an $e$-dimensional
Brownian motion with covariance matrix $\Sigma=\Cov(W(1))=\lim_{n\to
\infty}\Cov_\mu(W_n(1))$ given by
\[
\Sigma^{\beta\gamma}  =\int_\Lambda v^\beta
v^\gamma \,d\mu+ \sum_{n=1}^\infty
\int_\Lambda\bigl(v^\beta v^\gamma\circ
f^n+ v^\gamma v^\beta\circ f^n\bigr) \,d
\mu.
\]
\item[(ii)]
$\BBW^{\beta\gamma}(t)=\int_0^t W^\beta \,dW^\gamma+ E^{\beta\gamma
}t$ where
$E=\lim_{n\to\infty} \E_\mu(\BBW_n(1))$ is given by
\[
E^{\beta\gamma}= \sum_{n=1}^\infty\int
_\Lambda v^\beta v^\gamma \circ
f^n \,d\mu.
\]
\end{longlist}
\end{thmm}

Given $a\dvtx\R^d\to\R^d$, $b\dvtx\R^d\to\R^{d\times e}$, we define
$X_n\in D([0,\infty),\R^d)$, to be the solution to an appropriately
discretized version of equation \eqref{eq-ODE}.
Namely, we set
$X_n(t)=X_{[nt],n}$ where
\[
X_{j+1,n}=X_{j,n}+n^{-1}a(X_{j,n})+b(X_{j,n})
\biggl(W_n\biggl({\frac
{j+1}{n}}\biggr)-W_n\biggl({
\frac{j}{n}}\biggr) \biggr),\qquad X_{0,n}=\xi.
\]

\begin{thmm}[(Convergence to SDE, discrete time)] \label{thmm-approx2}
Suppose that $\Lambda\subset M$ is a mixing hyperbolic basic set with
equilibrium measure $\mu$, and that $v\dvtx\Lambda\to\R^e$ is H\"
older with
$\int_\Lambda v \,d\mu=0$.
Let $W_n$, $W$ and $E$ be as in Theorem~\ref{thmm-WW2}.
Let $a\dvtx\R^d\to\R^d$ be $C^{1+}$ and $b\dvtx\R^d\to\R
^{d\times e}$ be
$C^{2+}$, and define
$X_n\in D([0,\infty),\R^d)$ as above.

Then $X_n\to_w X$ in $D([0,\infty),\R^d)$ as $n\to\infty$, where $X$
satisfies the SDE
\[
dX= \biggl\{a(X)+\sum_{\alpha,\beta,\gamma}E^{\beta\gamma
}
\partial^\alpha b^{\beta}(X) b^{ \alpha\gamma}(X) \biggr\} \,dt+b(X)
\,dW, \qquad X(0)=\xi.
\]
\end{thmm}

\section{Cohomological invariance for iterated integrals}
\label{sec-cohom}

In this section, we present a result which is of independent
theoretical interest but which in particular significantly simplifies
the subsequent calculations.

Let $f\dvtx\Lambda\to\Lambda$ be an invertible or noninvertible map with
invariant probability measure $\mu$.
Suppose that
$v,\hat v\dvtx\Lambda\to\R^e$ are mean zero observables lying in $L^2$.
Define $W_n\in D([0,\infty),\R^e)$ and $\BBW_n\in D([0,\infty),\R
^{e\times e})$ as in \eqref{eq-WWW}, and similarly define
$\widehat W_n\in D([0,\infty),\R^e)$ and $\widehat\BBW_n\in
D([0,\infty
),\R^{e\times e})$ starting from $\hat v$ instead of~$v$.

We say that $v$ and $\hat v$ are \emph{$L^2$-cohomologous} if there exists
$\chi\dvtx\Lambda\to\R^e$ lying in $L^2$ such that
$v=\hat v+\chi\circ f-\chi$.
It is then easy to see that $W_n$ satisfies the WIP if and only if
$\widehat W_n$ satisfies the WIP and moreover the weak limits of
$W_n$ and $\widehat W_n$ coincide.
However, the weak limits of $\BBW_n$ and $\widehat\BBW_n$ need not coincide.
The following result supplies the correction factor needed to recover
identical weak limits.

\begin{thmm} \label{thmm-cohom}
Suppose that $f\dvtx\Lambda\to\Lambda$ is mixing and
that $v, \hat v\in L^2(\Lambda,\R^e)$ are $L^2$-cohomologous mean
zero observables.
Let $1\le\beta,\gamma\le e$.
Then the limit
$\lim_{n\to\infty}
\sum_{j=1}^n \int_\Lambda
(v^\beta v^\gamma\circ f^j-
\hat v^\beta\hat v^\gamma\circ f^j) \,d\mu$ exists and
\[
\BBW_n^{\beta\gamma}(t)- \widehat\BBW_n^{\beta\gamma}(t)
\to t \sum_{j=1}^\infty\int
_\Lambda\bigl(v^\beta v^\gamma\circ
f^j- \hat v^\beta\hat v^\gamma\circ
f^j\bigr) \,d\mu\qquad \mbox{a.e.},
\]
as $n\to\infty$, uniformly on compact subsets of $[0,\infty)$.

In particular, the weak limits of the processes
\[
\BBW_n^{\beta\gamma}(t)-t \sum_{j=1}^n
\int_\Lambda v^\beta v^\gamma\circ
f^j \,d\mu,\qquad \widehat\BBW_n^{\beta\gamma}(t)-t \sum
_{j=1}^n \int_\Lambda
\hat v^\beta\hat v^\gamma\circ f^j \,d\mu,
\]
coincide (in the sense that if one limit exists, then so does the other
and they are equal).
\end{thmm}

\begin{pf}
Write $v=\hat v+a$,
$a=\chi\circ f-\chi$, and $A_n(t)=n^{-1/2}\sum_{j=0}^{[nt]-1}a\circ f^j$.
Then
\[
\BBW_n^{\beta\gamma}(t)- \widehat\BBW_n^{\beta\gamma}(t)=
\int_0^t W_n^\beta
\,dW_n^\gamma -\int_0^t
\widehat W_n^\beta \,d\widehat W_n^\gamma
=\int_0^t A_n^\beta
\,dW_n^\gamma +\int_0^t
\widehat W_n^\beta \,dA_n^\gamma.
\]
Now
\begin{eqnarray*}
\int_0^t A_n^\beta
\,dW_n^\gamma & = &n^{-1}\sum
_{j=0}^{[nt]-1}\sum_{i=0}^{j-1}a^\beta
\circ f^i v^\gamma \circ f^j =n^{-1}
\sum_{j=0}^{[nt]-1}\bigl(\chi^\beta
\circ f^j-\chi^\beta\bigr) v^\gamma \circ
f^j
\\
& =& n^{-1}\sum_{j=0}^{[nt]-1}\bigl(
\chi^\beta v^\gamma\bigr)\circ f^j-
n^{-1}\chi^\beta\sum_{j=0}^{[nt]-1}v^\gamma
\circ f^j,
\end{eqnarray*}
which converges to $t\int_\Lambda\chi^\beta v^\gamma \,d\mu$ a.e. by
the ergodic theorem.

A similar argument for the remaining term, after changing order of
summation yields that $\int_0^t \widehat W_n^\beta \,dA_n^\gamma\to
-t\int_\Lambda\hat v^\beta\chi^\gamma\circ f \,d\mu$ a.e.

Hence, we have shown that
%
\begin{equation}
\label{eq-1} \BBW_n^{\beta\gamma}(t)- \widehat
\BBW_n^{\beta\gamma}(t)\to t \biggl(\int_\Lambda
\chi^\beta v^\gamma \,d\mu -\int_\Lambda\hat
v^\beta\chi^\gamma\circ f \,d\mu \biggr).
\end{equation}

Next,
\[
v^\beta v^\gamma\circ f^j-\hat v^\beta
\hat v^\gamma\circ f^j =\bigl(\chi^\beta\circ f-
\chi^\beta\bigr) v^\gamma\circ f^j + \hat
v^\beta\bigl(\chi^\gamma\circ f-\chi^\gamma\bigr)
\circ f^j,
\]
and so
%
\begin{eqnarray}
\label{eq-2}
\nonumber
&& \sum_{j=1}^n
\int_\Lambda v^\beta v^\gamma\circ
f^j \,d\mu -\sum_{j=1}^n \int
_\Lambda\hat v^\beta\hat v^\gamma\circ
f^j \,d\mu
\\
&&\qquad =\sum_{j=1}^n \int
_\Lambda \bigl\{ \bigl(\chi^\beta\circ f-
\chi^\beta\bigr) v^\gamma\circ f^j +\hat
v^\beta\bigl(\chi^\gamma\circ f-\chi^\gamma\bigr)
\circ f^j \bigr\}\,d\mu
\\
&&\qquad =\sum_{j=1}^n \int
_\Lambda \bigl\{ \bigl(\chi^\beta\circ
f^{n-j+1}-\chi^\beta\circ f^{n-j}\bigr)
v^\gamma\circ f^n \nonumber\\
&&\hspace*{98pt}{}+\hat v^\beta\bigl(
\chi^\gamma\circ f^{j+1}-\chi^\gamma\circ
f^j\bigr) \bigr\}\,d\mu\nonumber
\\
& &\qquad= \int_\Lambda\chi^\beta v^\gamma \,d\mu -
\int_\Lambda\hat v^\beta\chi^\gamma\circ f d
\mu+L_n,\nonumber
\end{eqnarray}
where
$L_n=\int_\Lambda(\hat v^\beta \chi^\gamma\circ f^{n+1}
-\chi^\beta v^\gamma\circ f^n) \,d\mu\to0$ as $n\to\infty$
by the mixing assumption.
The result is immediate from \eqref{eq-1} and \eqref{eq-2}.
\end{pf}

\begin{cor} \label{cor-cohom}
Let $f\dvtx\Lambda\to\Lambda$ be mixing and
let $v, \hat v\in L^2(\Lambda,\R^e)$ be $L^2$-cohomologous mean
zero observables.

Suppose that $(\widehat W_n,\widehat\BBW_n)\to_w (\widehat
W,\widehat
\BBW)$
in
$D([0,\infty),\R^e\times\R^{e\times e})$ as $n\to\infty$. Then
$(W_n,\BBW_n)\to_w (W,\BBW)$
in
$D([0,\infty),\R^e\times\R^{e\times e})$ as $n\to\infty$,
where $W=\widehat W$ and
\[
\BBW^{\beta\gamma}(t)= \widehat\BBW^{\beta\gamma}(t) +t \sum
_{j=1}^\infty\int_\Lambda
\bigl(v^\beta v^\gamma\circ f^j \,d\mu-
\hat{v}^\beta\hat v^\gamma\circ f^j\bigr) \,d\mu.
\]
\end{cor}

\begin{rmk}
For completeness, we describe the analogous result for semiflows. Again
the result is of independent theoretical interest even though we make
no use of it in this paper.

Let $\phi_t\dvtx\Omega\to\Omega$ be a mixing (semi)flow with invariant probability
measure $\nu$.
Suppose that
$v,\hat v\dvtx\Omega\to\R^e$ are mean zero observables lying in $L^2$.
Define $W_n$ and $\BBW_n$ as in \eqref{eq-W} and \eqref{eq-WW},
and similarly define
$\widehat W_n$ and $\widehat\BBW_n$ starting from $\hat v$ instead of $v$.

We say that $v$ and $\hat v$ are \emph{$L^2$-cohomologous} if there exists
$\chi\dvtx\Omega\to\R^e$ lying in $L^2$ such that\vspace*{1pt}
$\int_0^t v\circ\phi_s \,ds=\int_0^t\hat v\circ\phi_s \,ds+\chi\circ
\phi_t-\chi$.
Again, $W_n$ satisfies the WIP if and only if $\widehat W_n$ satisfies\vspace*{1pt}
the WIP and the weak limits coincide.
As in Theorem~\ref{thmm-cohom}, we find that
the limit
$\lim_{n\to\infty}\int_0^n \int_\Omega(v^\beta v^\gamma\circ
\phi_s-
\hat v^\beta\hat v^\gamma\circ\phi_s) \,d\nu$ exists and
\[
\BBW_n^{\beta\gamma}(t)- \widehat\BBW_n^{\beta\gamma}(t)
\to t \int_0^\infty\int_\Omega
\bigl(v^\beta v^\gamma\circ\phi_s- \hat
v^\beta\hat v^\gamma\circ\phi_s\bigr)\,d\nu \,ds\qquad
\mbox{a.e.},
\]
as $n\to\infty$, uniformly on compact subsets of $[0,\infty)$.
The proof is almost identical to that of Theorem~\ref{thmm-cohom}, and
hence is omitted.
\end{rmk}

\section{Iterated WIP for noninvertible maps}
\label{sec-NUE}

A sufficient condition for Theorem~\ref{thmm-WW2}
is that $f\dvtx\Lambda\to\Lambda$ is a mixing uniformly expanding map.
More generally, in this section we consider a class of nonuniformly
expanding maps with sufficiently rapid decay of correlations.
The underlying hypotheses can be satisfied only by noninvertible maps;
see Section~\ref{sec-NUH} for more general hypotheses appropriate for
invertible maps.

In Section~\ref{sec-mild} we give more details on the class of maps
that is considered in this section. In Section~\ref{sec-WWmap}, we
prove the iterated WIP for these maps.

\subsection{Noninvertible maps}
\label{sec-mild}

Let $f\dvtx\Lambda\to\Lambda$ be an ergodic measure-preserving map defined
on a probability space $(\Lambda,\mu)$ and let $v\dvtx\Lambda\to\R
^d$ be an
integrable observable with $\int_\Lambda v \,d\mu=0$.
Let $P\dvtx L^1(\Lambda)\to L^1(\Lambda)$ be the transfer operator
for $f$
given by
$\int_\Lambda Pw_1 w_2 \,d\mu=\int_\Lambda w_1 Uw_2 \,d\mu$ for
$w_1\in L^1(\Lambda)$,
$w_2\in L^\infty(\Lambda)$ where $Uw=w\circ f$.

\begin{defn} \label{def-decomp}
Let $p\ge1$.
We say that $v$ \emph{admits an $L^p$ martingale-coboundary decomposition}
if there exists $m,\chi\in L^p(\Lambda,\R^e)$ such that
%
\begin{equation}
\label{eq-decomp} v=m+\chi\circ f-\chi,\qquad m\in\ker P.
\end{equation}
We refer to $m$ as the \emph{martingale part} of the decomposition.
\end{defn}

\begin{rmk} \label{rmk-decomp}
The reason for calling $m$ a martingale will become clearer in
Section~\ref{sec-WWmap}.
For the time being, we note that it is standard and elementary that
$PU=I$ and
$UP=E( \cdot|f^{-1}\mathcal{B})$
where $\mathcal{B}$ is the underlying $\sigma$-algebra.
In particular
$E(m|f^{-1}\mathcal{B})=0$.
\end{rmk}

Our main result in this section is the following.

\begin{thmm} \label{thmm-WWmap} Suppose that $f$ is mixing and
that the decomposition \eqref{eq-decomp} holds with $p=2$.
Then the conclusion of Theorem~\ref{thmm-WW2} is valid.
\end{thmm}

\begin{prop} \label{prop-decomp}
Let $p\ge1$.
A sufficient condition for \eqref{eq-decomp} to hold is that
$v\in L^\infty$ and
there are constants $C>0$, $\tau>p$ such that
%
\begin{equation}
\label{eq-decay} \biggl|\int_\Lambda v w\circ f^n \,d\mu \biggr|
\le C\|w\|_\infty n^{-\tau}\qquad \mbox{for all } w\in
L^\infty, n\ge1.
\end{equation}
\end{prop}

\begin{pf}
By duality, $\|P^nv\|_1\le Cn^{-\tau}$.
Also, $\|P^nv\|_\infty\le\|v\|_\infty$ and it follows that
$\|P^nv\|_p\le\|v\|_\infty^{1-1/p}(Cn^{-\tau})^{1/p}$ which is summable.

Define $\chi=\sum_{n=1}^\infty P^nv\in L^p$, and write
$v=m +\chi\circ f-\chi$ where $m\in L^p$.
Applying $P$ to both sides and using the fact that $PU=I$,
we obtain that $m\in\ker P$.
\end{pf}

There are large classes of noninvertible maps for which the decay
condition~\eqref{eq-decay} has been established for sufficiently
regular $v$;
see Section~\ref{sec-gen}.
In particular, for uniformly expanding maps the decay is exponential
for H\"older continuous~$v$,
so $\tau$ and $p$ can be chosen arbitrarily large.

In the remainder of this subsection, we reduce Theorem~\ref{thmm-WWmap}
to the martingale part.
Define the cadlag processes $M_n\in D([0,\infty),\R^e)$,
$\BBM_n\in D([0,\infty),\break \R^{e\times e})$,
\begin{eqnarray*}
M_n(t) & =&n^{-1/2}\sum_{j=0}^{[nt]-1}m
\circ f^j,
\\
\BBM^{\beta\gamma}_n(t) & =&\int_0^t
M_n^\beta \,dM_n^\gamma=n^{-1}
\sum_{0\le i<j\le[nt]-1} m^\beta\circ f^i
m^\gamma\circ f^j.
\end{eqnarray*}
Theorem~\ref{thmm-WWmap} follows from the following lemma.

\begin{lemma} \label{lem-WWmap}
Suppose that $f$ is ergodic and that $m\in L^2(\Lambda,\R^e)$ with
$Pm=0$. Then
$(M_n,\BBM_n)\to_w (W,I)$ in $D([0,\infty),\R^e\times\R^{e\times e})$,
as $n\to\infty$, where $W$ is an $e$-dimensional Brownian motion with
covariance matrix\break $\Cov(W(1))=\int_\Lambda m m^T \,d\mu$ and
$I^{\beta\gamma}(t)=\int_0^t W^\beta \,dW^\gamma$.
\end{lemma}

\begin{pf*}{Proof of Theorem~\ref{thmm-WWmap}}
We apply Corollary~\ref{cor-cohom} with $\hat v=m$.
Note that $\int_\Lambda m m^T\circ f^j \,d\mu=
\int_\Lambda P^jm m^T \,d\mu=0$ for all $j\ge1$.
By Theorem~\ref{thmm-cohom}, $E=\sum_{j=1}^\infty v v^T\circ f^j \,d\mu$
is a convergent series.
By Corollary~\ref{cor-cohom}, $(W_n,\BBW_n)\to_w(W,\BBW)$ where
$\Sigma=\Cov(W(1))=\int_\Lambda m m^T \,d\mu$ and $\BBW(t)=I(t)+Et$.

It remains to prove that $\Sigma^{\beta\gamma}=\lim_{n\to\infty
}\Cov
^{\beta\gamma}_\mu(W_n(1))=
\int_\Lambda v^\beta v^\gamma \,d\mu+
\sum_{n=1}^\infty\int_\Lambda(v^\beta v^\gamma\circ f^n+
v^\gamma v^\beta\circ f^n) \,d\mu$ and
that $E=\lim_{n\to\infty} \E_\mu(\BBW_n(1))$.

Define $v_n=\sum_{j=0}^{n-1}v\circ f^j$,
$m_n=\sum_{j=0}^{n-1}m\circ f^j$.
Then
\[
\int_\Lambda m_nm_n^T \,d
\mu= \sum_{0\le i,j\le n-1} \int_\Lambda m\circ
f^i\bigl(m\circ f^j\bigr)^T \,d\mu=n\Sigma.
\]
Equivalently, $c^T\Sigma c=n^{-1}\int_\Lambda(c^Tm_n)^2 \,d\mu$
for all $c\in\R^e$, $n\ge1$.
Let $\|\cdot \|_2$ denote the $L^2$ norm on $(\Lambda,\mu)$. We have that
$n^{{1}/2}(c^T\Sigma c)^{{1}/2}= \|c^Tm_n\|_2$. By \eqref{eq-decomp},
$v_n-m_n=\chi\circ f^n-\chi$. Using $f$-invariance of
$\mu$,
\begin{eqnarray*}
\bigl|\bigl\|c^Tv_n\bigr\|_2-n^{{1}/2}
\bigl(c^T\Sigma c\bigr)^{{1}/2} \bigr|& =&\bigl |\bigl\|c^Tv_n
\bigr\|_2-\bigl\|c^Tm_n\bigr\|_2 \bigr|\le
\bigl\|c^T(v_n-m_n)\bigr\|_2\\
&\le& 2
\bigl\|c^T\chi\bigr\|_2,
\end{eqnarray*}
and hence $\lim_{n\to\infty}n^{-{1}/2}\|c^Tv_n\|_2=(c^T\Sigma
c)^{{1}/2}$.
Equivalently,
%
\begin{equation}
\label{eq-Sigma1} \Sigma=\lim_{n\to\infty}n^{-1}\int
_\Lambda v_n v_n^T \,d\mu=
\lim_{n\to
\infty}\Cov_\mu\bigl(W_n(1)\bigr).
\end{equation}

Let $a_r=\int_\Lambda v\circ f^r v^T \,d\mu$ and $s_k=\sum_{r=1}^k a_r$.
Compute that
\begin{eqnarray*}
\sum_{0\le j<i\le n-1}\int_\Lambda v\circ
f^{i-j} v^T \,d\mu & =& \sum_{1\le r<n}(n-r)
\int_\Lambda v\circ f^r v^T \,d\mu
\\
& = &\sum_{1\le r<n}(n-r)a_r =\sum
_{k=1}^n s_k.
\end{eqnarray*}
Hence,
%
\begin{eqnarray}
\label{eq-Sigma2}
\lim_{n\to\infty}n^{-1}\sum
_{0\le j<i\le n-1}\int_\Lambda v\circ
f^{i-j} v^T \,d\mu & =&\lim_{n\to\infty}n^{-1}
\sum_{k=1}^n s_k =\lim
_{n\to\infty}s_n
\nonumber
\\[-8pt]
\\[-8pt]
\nonumber
& =&\sum_{r=1}^\infty\int
_\Lambda v\circ f^r v^T \,d\mu.
\end{eqnarray}
Similarly,
%
\begin{equation}
\label{eq-Sigma3} \lim_{n\to\infty}n^{-1}\sum
_{0\le i<j\le n-1}\int_\Lambda v \bigl(v\circ
f^{j-i}\bigr)^T \,d\mu = \sum_{r=1}^\infty
\int_\Lambda v \bigl(v\circ f^r
\bigr)^T \,d\mu.
\end{equation}

Write
\begin{eqnarray*}
n^{-1}\int_\Lambda v_nv_n^T
\,d\mu & =& n^{-1} \sum_{0\le i,j\le n-1}\int
_\Lambda v\circ f^i\bigl(v\circ f^j
\bigr)^T \,d\mu
\\
& =& \int_\Lambda vv^T \,d\mu + n^{-1}
\sum_{0\le j<i\le n-1}\int_\Lambda v\circ
f^{i-j} v^T \,d\mu
\\
&&{} + n^{-1} \sum_{0\le i<j\le n-1}\int
_\Lambda v \bigl(v\circ f^{j-i}\bigr)^T \,d
\mu.
\end{eqnarray*}
By \eqref{eq-Sigma1}, \eqref{eq-Sigma2}, \eqref{eq-Sigma3},
$\Sigma=\int_\Lambda v v^T \,d\mu+
\sum_{r=1}^\infty\int_\Lambda(v\circ f^r v^T+
v (v\circ f^r)^T) \,d\mu$.

Finally,
$\E_\mu(\BBW_n(1))=n^{-1}\sum_{0\le i<j\le n-1}\int_\Lambda v
(v\circ
f^{j-i})^T \,d\mu$, so it follows from~\eqref{eq-Sigma3} that
$\lim_{n\to\infty}\E_\mu(\BBW_n(1))=E$.
\end{pf*}

\subsection{Proof of Lemma \texorpdfstring{\protect\ref{lem-WWmap}}{4.5}}
\label{sec-WWmap}

\begin{rmk} \label{rmk-WIP}
The $M_n\to_w W$ part of Lemma~\ref{lem-WWmap} is standard but we give
the proof for completeness.
The statement can be obtained from the proof of Lemma~\ref{lem-WWmap}
by ignoring the $\BBM_n$ component. In particular, our use of this fact
in the proof of Lemma~\ref{lem-gh} below is not circular.
\end{rmk}

Recall that $m$ is $\mathcal{B}$-measurable and
$m\in\ker P$ so $E(m|f^{-1}\mathcal{B})=0$.
Similarly, $m\circ f^j$ is $f^{-j}\mathcal{B}$-measurable and
$E(m\circ f^j|f^{-(j+1)}\mathcal{B})
=E(m|f^{-1}\mathcal{B})\circ\break f^j=0$.
If the sequence of $\sigma$-algebras $f^{-j}\mathcal{B}$
formed a filtration, then $M_n$ would be a martingale and we could
apply Kurtz and Protter \cite{KurtzProtter91}, Theorem~2.2
(see also \cite{JakubowskiMeminPages89})
to obtain a limit for $(M_n,\BBM_n)$.

In fact, the $\sigma$-algebras are decreasing: $f^{-j}\mathcal
{B}\supset f^{-(j+1)}\mathcal{B}$ for all $j$.
To remedy this, we pass to the natural extension $\tilde f\dvtx\tilde
\Lambda
\to\tilde\Lambda$. This is an invertible map with ergodic invariant
measure $\tilde\mu$, and there is a measurable projection $\pi\dvtx
\tilde
\Lambda\to\Lambda$ such that
$\pi\tilde f=f\pi$ and $\pi_*\tilde\mu=\mu$.
The observable $m\dvtx\Lambda\to\R^e$ lifts to an observable
$\tilde m=m\circ\pi\dvtx\tilde\Lambda\to\R^e$ and the joint distributions
of $\{m\circ f^j\dvtx j\ge0\}$ are identical to those of
$\{\tilde m\circ\tilde f^j\dvtx j\ge0\}$.

Define
\begin{eqnarray*}
\widetilde {M}_n(t)&=&n^{-1/2}\sum_{j=0}^{[nt]-1}
\tilde m\circ\tilde f^j,
\\
\widetilde{\mathbb{M}}^{\beta\gamma}_n(t)&=&
\int_0^t \widetilde M^\beta_n
\,d\widetilde M^\gamma_n =n^{-1} \sum
_{0\le i<j\le[nt]-1} \tilde m^\beta\circ\tilde f^i
\tilde m^\gamma\circ\tilde f^j.
\end{eqnarray*}
Then $(\widetilde M_n,\widetilde\BBM_n)=(M_n,\BBM_n)\circ\pi$ and
$\pi
$ is measure preserving, so
it is equivalent to prove that
%
\begin{equation}
\label{eq-back} (\widetilde M_n,\widetilde\BBM_n)
\to_w (W,I)\qquad \mbox{in }\mbox D\bigl([0,\infty),\R^e
\times\R^{e\times e}\bigr).
\end{equation}

Let $\widetilde{\mathcal{B}}=\pi^{-1}\mathcal{B}$.
Again $\tilde f^{-j}\widetilde{\mathcal{B}}\supset
\tilde f^{-(j+1)}\widetilde{\mathcal{B}}$ but this means that
$\{\mathcal{F}_j, {j\ge1}\}=\{\tilde f^j\widetilde{\mathcal{B}},
j\ge
1\}$ is an \emph{increasing} sequence of
$\sigma$-algebras.
Moreover, $\tilde m\circ\tilde f^{-j}$ is $\mathcal{F}_j$-measurable
and $E(\tilde m\circ\tilde f^{-j}|\mathcal{F}_{j-1})=0$.
Hence, the ``backward'' process
\[
\widetilde M_n^-(t)=n^{-{1}/2}\sum
_{j=-[nt]}^{-1}\tilde m\circ \tilde f^j
\]
forms an ergodic stationary martingale. Similarly, define
\[
\widetilde\BBM_n^{\beta\gamma,-}(t)=\int_0^t
\widetilde M_n^{\beta,-} \,d\widetilde M_n^{\gamma,-}=
n^{-1} \sum_{[-nt]\le j<i\le-1} \tilde m^\beta
\circ\tilde f^i \tilde m\circ\tilde f^j.
\]
Note that $\int_{\tilde\Lambda} \tilde m \tilde m^T \,d\tilde\mu
=\int_\Lambda m m^T \,d\mu$.

\begin{prop} \label{prop-KP}
$(\widetilde M_n^-,\widetilde\BBM_n^-)\to_w (W,I) \mbox{ in }
D([0,\infty),\R^e\times\R^{e\times e})$ as\break \mbox{$n\to\infty$}.
\end{prop}

\begin{pf} We verify the hypotheses of
Kurtz and Protter \cite{KurtzProtter91}, Theorem~2.2
(with \mbox{$\delta=\infty$} and $A_n\equiv0$).
We have already seen that
$\widetilde M_n^-$ is a martingale.
Also, by the calculation in the proof of
Theorem~\ref{thmm-WWmap}, $E(\widetilde M_n^{\gamma,-}(t)^2)=n^{-1}\|
\sum_{j=1}^{[nt]}\tilde m^\gamma\circ\tilde f^{-j}\|_2^2
=t\int_{\tilde\Lambda} (\tilde m^\gamma)^2 \,d\tilde\mu$
independent of
$n$, so condition C2.2(i)
in \cite{KurtzProtter91}, Theorem~2.2, is trivially satisfied.

The WIP for stationary ergodic $L^2$ martingales
(e.g., \cite{Brown71,McLeish74}) implies that $\widetilde M_n^-\to_w W$
in $D([0,\infty),\R^e)$.
In particular,
$(\widetilde M_n^{\beta,-},\widetilde M_n^{\gamma,-})\to_w (W^\beta
,W^\gamma)$ in $D([0,\infty),\R^2)$.
Hence, the result follows from \cite{KurtzProtter91}, Theorem~2.2.
\end{pf}

It remains to relate weak convergence of $(\widetilde M_n^-,\widetilde
\BBM_n^-)$
and $(\widetilde M_n,\widetilde\BBM_n)$. It suffices to work in
$D([0,T],\R^e\times\R^{e\times e})$
for each fixed integer $T\ge1$.

\begin{lemma} \label{lem-gh}
Let $g(u)(t)=u(T)-u(T-t)$
and $h(u,v)(t)=u(T-t)(v(T)-v(T-t))$.
Let $\ast$ denote matrix transpose in $\R^{e\times e}$.
Then
\[
(\widetilde M_n,\widetilde\BBM_n)\circ\tilde
f^{-nT}= \bigl(g\bigl(\widetilde M_n^-\bigr), \bigl(g\bigl(
\widetilde\BBM_n^-\bigr)-h\bigl(\widetilde M_n^-\bigr)
\bigr)^* \bigr)+F_n,
\]
where $\sup_{t\in[0,T]}F_n(t)\to0$ a.e.
\end{lemma}

\begin{pf}
In this proof, we suppress the tildes.
First, we show that
$M_n\circ f^{-nT}=g(M_n^-)+F^0_n$,
where $\sup_{t\in[0,T]}F^0_n(t)\to0$ a.e.

We have
\begin{eqnarray*}
M_n(t)\circ f^{-nT} & =&n^{-{1}/2} \sum
_{j=0}^{[nt]-1}m\circ f^j\circ
f^{-nT} =n^{-{1}/2} \sum_{j=-nT}^{[nt]-1-nT}m
\circ f^j
\\
& =& M_n^-(T)-M_n^-(T-t)+F^0_n(t).
\end{eqnarray*}
Here, $F^0_n$ consists of at most one term and we can write
\[
\bigl|F^0_n(t)\bigr|\le n^{-{1}/2}\Bigl |\max
_{j=1,\ldots,nT} m\circ f^{-j} \Bigr|.
\]
It suffices to work componentwise, so suppose without loss that $e=1$.
By the ergodic theorem, $n^{-1}\sum_{j=1}^n m^2\circ f^{-j}\to\int_{\Lambda} m^2 \,d\mu$, and
so $n^{-1}m^2\circ f^{-n}\to0$.
It follows that $n^{-1}\max_{j=0,\ldots,nT}m^2\circ f^{-j}\to0$ a.e.
and so $\sup_{t\in[0,T]}F^0_n(t)\to0$ a.e.

Next, we show that
$\BBM_n\circ f^{-nT}=
 (g(\BBM_n^-)-h(M_n^-) )^*+F_n$,
where\break \mbox{$\sup_{t\in[0,T]}F_n(t)\to0$} a.e.
We have
\begin{eqnarray*}
\BBM_n^{\beta\gamma}(t) & =&n^{-1}\sum
_{j=0}^{[nt]-1} \Biggl(\sum_{i=0}^{j-1}m^\beta
\circ f^i \Biggr) m^\gamma\circ f^j,
\\
\BBM_n^{\beta\gamma,-}(t) & =&n^{-1}\sum
_{j=-[nt]+1}^{-1} \Biggl(\sum_{i=[-nt]}^{j-1}m^\gamma
\circ f^i \Biggr) m^\beta\circ f^j.
\end{eqnarray*}
Hence,
%
\begin{eqnarray}
\label{eq-IEF} &&\hspace*{-4pt}\BBM_n^{\beta\gamma}(t)\circ f^{-nT}\nonumber \\
&&\hspace*{-7pt}\qquad=
n^{-1}\sum_{j=-nT}^{[nt]-1-nT} \sum
_{i=-nT}^{j-1} m^\beta\circ
f^i m^\gamma\circ f^j
\nonumber
\\[-8pt]
\\[-8pt]
\nonumber
&&\hspace*{-7pt}\qquad = n^{-1} \Biggl(\sum_{j=-nT}^{-nT}
+\sum_{j=-nT+1}^{-1} - \sum
_{j=[nt]-nT+1}^{-1} - \sum_{j=[nt]-nT}^{[nt]-nT}
\Biggr) \sum_{i=-nT}^{j-1} m^\beta
\circ f^i m^\gamma\circ f^j
\\
&&\hspace*{-7pt}\qquad = F^1_n(t)+\BBM_n^{\gamma\beta,-}(T)-E_n(t)-F^2_n(t),\nonumber
\end{eqnarray}
where
\begin{eqnarray*}
F^1_n(t) & =&n^{-1}\sum
_{i=-nT}^{-nT-1}m^\beta\circ f^i
m^\gamma\circ f^{-nT},
\\
F^2_n(t) & = &\Biggl(n^{-{1}/2} \sum
_{i=-nT}^{[nt]-nT-1} m^\beta\circ f^i
\Biggr) \bigl(n ^{-{1}/2}m^\gamma\circ f^{[nt]-nT} \bigr),
\\
E_n(t) & = &n^{-1}\sum_{j=[nt]-nT+1}^{-1}
\sum_{i=-nT}^{j-1} m^\beta\circ
f^i m^\gamma\circ f^j.
\end{eqnarray*}
Now $F^1_n(t)$ consists of only two terms and clearly converges to $0$
almost everywhere.
The first factor in $F^2_n$ converges weakly to $W^\beta$ (see
Remark~\ref{rmk-WIP}) and the second factor converges to $0$ almost
everywhere by the ergodic theorem. Hence,
$\sup_{t\in[0,T]}Z|F^r_n(t)|\to0$ a.e. for $r=1,2$.
Moreover,
%
\begin{eqnarray}
\label{eq-EH}
E_n(t) & =& n^{-1}\sum
_{j=[nt]-nT+1}^{-1} \Biggl(\sum_{i=-nT}^{-nT+[nt]-1}
+\sum_{i=-nT+[nt]}^{j-1} \Biggr) m^\beta
\circ f^i m^\gamma\circ f^j
\nonumber
\\[-8pt]
\\[-8pt]
\nonumber
& = &H_n(t)+\BBM_n^{\gamma\beta,-}(T-t)+F^3_n(t),
\end{eqnarray}
where
%
\begin{eqnarray}
\label{eq-H}
H_n(t) & =& \Biggl( n^{-{1}/2} \sum
_{j=[nt]-nT}^{-1} m^\gamma\circ
f^j \Biggr) \Biggl( n^{-{1}/2} \sum
_{i=-nT}^{-nT+[nt]-1} m^\beta\circ f^i
\Biggr)
\nonumber
\\[-8pt]
\\[-8pt]
\nonumber
& =& M_n^{\gamma,-}(T-t) \bigl(M_n^{\beta,-}(T)-M_n^{\beta
,-}(T-t)
\bigr),
\end{eqnarray}
and $F^3_n(t)=n^{-1}\sum_{i=-nT}^{-nT+[nt]-1}m^\beta\circ f^i
m^\gamma
\circ f^{[nt]-nT+1}$.
Again, $\sup_{t\in[0,T]}|F^3_n(t)|\to0$ a.e. by the ergodic theorem.
The result follows from \eqref{eq-IEF}, \eqref{eq-EH}, \eqref{eq-H}.
\end{pf}

\begin{prop} \label{prop-caglad}
Let
$\widetilde D([0,T],\R^q)$ denote the space of caglad functions from
$[0,T]$ to $\R^q$
with the standard Skorokhod $\mathcal{J}_1$ topology.
Suppose that $A_n=B_n+F_n$ where
$A_n\in D([0,T],\R^q)$,
$B_n\in\widetilde D([0,T],\R^q)$,
and $F_n\to0$ uniformly in probability.
If $Z$ has continuous sample paths and
$B_n\to_w Z$ in $\widetilde D([0,T],\R^q)$, then
$A_n\to_w Z$ in $D([0,T],\R^q)$.
\end{prop}

\begin{pf}
It is clear that the limiting finite distributions of $A_n$ coincide
with those
of $B_n$, so it suffices to show that $A_n$ inherits tightness from $B_n$.
One way to see this is to consider the following Arzela--Ascoli-type
characterization \cite{Skorokhod56}, valid in both $D([0,T],\R^q)$ and
$\widetilde D([0,T],\R^q)$.

Tightness of $B_n$ in $\widetilde D([0,T],\R^q)$ implies that
for any $\varepsilon>0$, $k\ge1$, there exists $C>0$, $\delta_k>0$,
$n_k\ge1$ such that $P(|B_n|_\infty>C) < \varepsilon$ for all $n\ge
1$ and
$P(\omega(B_n,\delta_k)>1/k) < \varepsilon$ for all $n\ge n_k$,
where
\[
\omega(\psi,\delta)=\sup_{t-\delta<t'<t<t''<t+\delta}\min\bigl\{\bigl|\psi (t)-\psi
\bigl(t'\bigr)\bigr|,\bigl|\psi(t)-\psi\bigl(t''
\bigr)\bigr|\bigr\}
\]
(where $t,t',t''$ are restricted to $[0,T]$).
These criteria are also satisfied by $F_n$ for trivial reasons, and hence
by $A_n$ establishing tightness of $A_n$ in $D([0,T],\R^q)$.
\end{pf}

\begin{cor} \label{cor-gh}
$(\widetilde M_n,\widetilde\BBM_n)\to_w  (g(W),
 (g(I)-h(W) )^* )$ in
$D([0,T], \R^e\times\R^{e\times e})$ as $n\to\infty$.
\end{cor}

\begin{pf}
Recalling the notation from Lemma~\ref{lem-gh},
observe that the functional $\chi\dvtx D([0,T],\R^e\times\R
^{e\times e})\to
\widetilde D([0,T],\R^e\times\R^{e\times e})$ given by
$\chi(u,v)=(g(u),\break (g(v)-h(u))^*)$ is continuous.
Hence, it follows from Proposition~\ref{prop-KP}\vspace*{1pt} and the continuous
mapping theorem that
$(g(\widetilde M_n^-),
(g(\widetilde\BBM_n^-)-h(\widetilde M_n^-))^*)
\to_w (g(W),\break (g(I)-h(W))^*)$ in $\widetilde D([0,T],\R^e\times\R
^{e\times e})$.
The result is now immediate from Lemma~\ref{lem-gh} and Proposition~\ref
{prop-caglad}.
\end{pf}

\begin{lemma} \label{lem-identify}
$(g(W),(g(I)-h(W))^*)=_d(W,I)$ in
$D([0,T],\R^e\times\R^{e\times e})$.
\end{lemma}

\begin{pf}
\textit{Step} 1. $g(W)=_d W$
in $D([0,T],\R^e)$.
To see this, note that both processes are
Gaussian with continuous sample
paths and $g(W)(0)=W(0)=0$. One easily verifies that $\Cov
(g(W)(t_1),g(W)(t_2))=t_1\Sigma$ for all $0\le t_1\le t_2\le T$.
Hence, $g(W)=_d W$.

\textit{Step} 2. Introduce the process $J(t)=\int_0^t g(W) \,dg(W)$.
We claim that\break $(g(W), J)=_d (W,I)$. To see this, let
$Y_n(t)=\sum_{j=0}^{[nt]-1}W(j/n)(W((j+1)/n)-W(j/n))$ so
$(W,Y_n)\to_w (W,I)$.
Similarly,\vspace*{1pt} let
$Z_n(t)=\sum_{j=0}^{[nt]-1}g(W)(j/n)\*(g(W)((j+1)/n)-g(W)(j/n))$ so
$(g(W),Z_n)\to_w (g(W),J)$.
It is clear that $(W,Y_n)=_d(g(W),Z_n)$ so the claim follows.

\textit{Step} 3. We complete the proof by showing that
$J=(g(I)-h(W))^*$. Let $1\le\beta,\gamma\le e$. We show that
$g(I)^{\beta\gamma}-h(W)^{\beta\gamma}=J^{\gamma\beta}$.

Now
$J^{\gamma\beta}(t) = \int_0^t g(W)^\gamma \,dg(W)^\beta
=\lim_{n\to\infty} S_n$ where the limit is in probability and
\begin{eqnarray*}
S_n & =& \sum_{k=0}^{[nt]-1}g(W)^\gamma
\biggl({\frac{k}{n}}\biggr) \biggl(g(W)^\beta\biggl({
\frac{k+1}{n}}\biggr)- g(W)^\beta\biggl({\frac{k}{n}}\biggr)
\biggr)
\\
& =&\sum_{k=0}^{[nt]-1} \biggl(W^\gamma(T)-W^\gamma
\biggl(T-{\frac
{k}{n}}\biggr) \biggr)\\
&&\hspace*{28pt}{}\times \biggl(W^\beta\biggl(T-{
\frac{k}{n}}\biggr)- W^\beta\biggl(T-{\frac
{k+1}{n}}\biggr)
\biggr)
\\
& =&\sum_{k=0}^{[nt]-1}\sum
_{j=0}^{k-1} \biggl(W^\gamma\biggl(T-{
\frac
{j}{n}}\biggr)-W^\gamma\biggl(T-{\frac{j+1}{n}}\biggr)
\biggr)\\
&&\hspace*{46pt}{}\times \biggl(W^\beta\biggl(T-{\frac{k}{n}}\biggr)-
W^\beta\biggl(T-{\frac
{k+1}{n}}\biggr) \biggr)
\\
& =&\sum_{j=0}^{[nt]-2}\sum
_{k=j+1}^{[nt]-1} \biggl(W^\beta\biggl(T-{
\frac{k}{n}}\biggr)- W^\beta\biggl(T-{\frac
{k+1}{n}}\biggr)
\biggr) \\
&&\hspace*{58pt}{}\times \biggl(W^\gamma\biggl(T-{\frac{j}{n}}\biggr)-W^\gamma
\biggl(T-{\frac
{j+1}{n}}\biggr) \biggr)
\\
& =&\sum_{j=0}^{[nt]-2} \biggl(W^\beta
\biggl(T-{\frac{j+1}{n}}\biggr)- W^\beta\biggl(T-{
\frac
{[nt]}{n}}\biggr) \biggr) \\
&&\hspace*{29pt}{}\times \biggl(W^\gamma\biggl(T-{
\frac{j}{n}}\biggr)-W^\gamma\biggl(T-{\frac
{j+1}{n}}\biggr)
\biggr).
\end{eqnarray*}

On the other hand,
$\{g(I)-h(W)\}^{\beta\gamma}(t) = \int_{T-t}^T (W^\beta-W^\beta(T-t))
\,dW^\gamma=\lim_{n\to\infty} T_n$ where
\begin{eqnarray*}
T_n & =&\sum_{i=[n(T-t)]}^{nT-1}
\biggl(W^\beta\biggl({\frac
{i}{n}}\biggr)-W^\beta(T-t)
\biggr) \biggl(W^\gamma\biggl({\frac{i+1}{n}}\biggr)-W^\gamma
\biggl({\frac{i}{n}}\biggr) \biggr)
\\
& =&\sum_{j=0}^{-[-nt]-1} \biggl(W^\beta
\biggl(T-{\frac{j+1}{n}}\biggr)-W^\beta (T-t) \biggr)\\
&&\hspace*{46pt}{}\times
\biggl(W^\gamma\biggl(T-{\frac{j}{n}}\biggr)-W^\gamma
\biggl(T-{\frac
{j+1}{n}}\biggr) \biggr).
\end{eqnarray*}

We claim that $\lim_{n\to\infty}(T_n-S_n)=0$ a.e. from which the
result follows.
When $nt$ is an integer, $S_n=T_n$. Otherwise,
$T_n-S_n= A_n+B_n$ where
\begin{eqnarray*}
A_n & =&\sum_{j=0}^{[nt]-2}
\biggl(W^\beta(T-t)-W^\beta\biggl(T-{\frac
{[nt]}{n}}
\biggr) \biggr) \biggl(W^\gamma\biggl(T-{\frac{j}{n}}
\biggr)-W^\gamma\biggl(T-{\frac
{j+1}{n}}\biggr) \biggr)
\\
& = &\biggl(W^\beta(T-t)-W^\beta\biggl(T-{\frac{[nt]}{n}}
\biggr) \biggr) \biggl(W^\gamma(T)-W^\gamma\biggl(T-\biggl({
\frac{[nt]-1}{n}}\biggr) \biggr)\biggr)
\end{eqnarray*}
and
\begin{eqnarray*}
&&B_n= \biggl(W^\beta\biggl(T-\biggl({\frac{[nt]+1}{n}}
\biggr)-W^\beta(T-t) \biggr)\biggr)\\
&&\hspace*{33pt}{}\times \biggl(W^\gamma\biggl(T-{
\frac{[nt]}{n}}\biggr)-W^\gamma\biggl(T-\biggl({\frac
{[nt]+1}{n}}
\biggr) \biggr)\biggr).
\end{eqnarray*}
The claim follows since $A_n\to0$ and $B_n\to0$ as $n\to\infty$.
\end{pf}

\begin{pf*}{Proof of Lemma~\ref{lem-WWmap}}
This follows from Corollary~\ref{cor-gh} and Lemma~\ref{lem-identify}.
\end{pf*}

\section{Iterated WIP for invertible maps}
\label{sec-NUH}

In this section, we prove an iterated WIP for invertible maps, and as a
special case we prove Theorem~\ref{thmm-WW2}.

For an invertible map $f\dvtx\Lambda\to\Lambda$, the transfer
operator $P$
is an isometry on $L^p$ for all $p$, so the hypotheses used in
Section~\ref{sec-NUE} are not applicable. We require the following more
general setting.

Suppose that in addition to the underlying probability space $(\Lambda
,\mu)$ and measure-preserving map $f\dvtx\Lambda\to\Lambda$, there
is an
additional probability space $(\bar\Lambda,\bar\mu)$ and
measure-preserving map $\bar f\dvtx\bar\Lambda\to\bar\Lambda$,
and there is
a semiconjugacy
$\pi\dvtx\Lambda\to\bar\Lambda$ with $\pi_*\mu=\bar\mu$ such that
$\pi\circ f=\bar f\circ\pi$. (The system on $\bar\Lambda$ is
called a
factor of
the system on $\Lambda$.)
We let $P$ denote the transfer operator for $\bar f$.

\begin{defn}
Let $v\dvtx\Lambda\to\R^e$ be of mean zero and
let $p\ge1$.
We say that $v$ \emph{admits an $L^p$ martingale-coboundary decomposition}
if there exists $m,\chi\in L^p(\Lambda,\R^e)$, $\bar m\in L^p(\bar
\Lambda,\R^e)$, such that
%
\begin{equation}
\label{eq-decomp_inv} v=m+\chi\circ f-\chi,\qquad m=\bar m\circ\pi,\qquad \bar m\in\ker P.
\end{equation}
\end{defn}

The definition is clearly more general than Definition~\ref
{def-decomp}, but
the consequences are unchanged.

\begin{thmm} \label{thmm-WWmap_inv} Suppose that $f$ is mixing and
that the decomposition \eqref{eq-decomp_inv} holds with $p=2$.
Then the conclusion of Theorem~\ref{thmm-WW2} is valid.
\end{thmm}

\begin{pf}
By Theorem~\ref{thmm-cohom}, we again reduce to considering the
martingale part~$m$.
Define the cadlag processes $(M_n,\BBM_n)$ and
$(\overline{M}_n,\overline{\BBM}_n)$ starting from
$m$ and $\bar m$, respectively.
Then $(M_n,\BBM_n)= (\overline{M}_n,\overline{\BBM}_n)\circ\pi$.
Hence, we reduce to
proving the iterated WIP for
$(M_n,\BBM_n)= (\overline{M}_n,\overline{\BBM}_n)$.
Since $\bar m\in\ker P$, we are now in the situation of Section~\ref
{sec-NUE}, and the result follows from Lemma~\ref{lem-WWmap}.
\end{pf}

For the remainder of this paper, hypotheses about the existence of a
martingale-coboundary decomposition refer only to the more general decomposition
in \eqref{eq-decomp_inv}.

\subsection{Applications of Theorem \texorpdfstring{\protect\ref{thmm-WWmap_inv}}{5.2}}

We consider first the case of
Axiom A (uniformly hyperbolic) diffeomorphisms.
By Bowen \cite{Bowen75}, any (nontrivial) hyperbolic basic set can be
modeled by a two-sided subshift of finite type $f\dvtx\Lambda\to
\Lambda$.
The alphabet consists of $k$ symbols $\{0,1,\ldots,k-1\}$ and there is a
transition matrix $A\in\R^{k\times k}$ consisting of zeros and ones.
The phase space $\Lambda$ consists of bi-infinite sequences
$y=(y_i)\in
\{0,1,\ldots,k-1\}^\Z$ such that $A_{y_i,y_{i+1}}=1$ for all $i\in\Z
$, and
$f$ is the shift $(fy)_i=y_{i+1}$.

For any $\theta\in(0,1)$, we define the metric $d_\theta(x,y)=\theta
^{s(x,y)}$
where the separation time $s(x,y)$ is the greatest integer $n\ge0$
such that
$x_i=y_i$ for $|i|\le n$.
Define $F_\theta(\Lambda)$ to be the space of $d_\theta$-Lipschitz functions
$v\dvtx\Lambda\to\R^e$ with Lipschitz constant $|v|_\theta=\sup_{x\neq
y}|x-y|/d_\theta(x,y)$ and norm $\|v\|_\theta=|v|_\infty+|v|_\theta$
where $|v|_\infty$ is the sup-norm. For each $\theta$, this norm makes
$F_\theta(\Lambda)$ into a Banach space.

As usual, we have the corresponding one-sided shift $\bar f\dvtx\bar
\Lambda
\to\bar\Lambda$
where $\bar\Lambda=\{0,1,\ldots,k-1\}^{\{0,1,2,\ldots\}}$, and the
associated function space
$F_\theta(\bar\Lambda)$. There is a natural projection $\pi\dvtx
\Lambda\to
\bar\Lambda$ that is a semiconjugacy between the shifts $f$ and $\bar
f$, and Lipschitz
observables $\bar v\in F_\theta(\bar\Lambda)$ lift to Lipschitz
observables $v=\bar v\circ\pi\in F_\theta(\Lambda)$.

A $k$-cylinder in $\bar\Lambda$ is a set of the form
$[a_0,\ldots,a_{k-1}]=\{y\in\bar\Lambda\dvtx y_i=a_i$ for all
$i=0,\ldots,k-1\}$, where $a_0,\ldots,a_{k-1}\in\{0,1,\ldots,k-1\}$.
The underlying $\sigma$-algebra $\overline{\mathcal{B}}$ is defined to
be the
$\sigma$-algebra generated by the $k$-cylinders. Note that $\bar{f}\dvtx
\bar
\Lambda\to\bar\Lambda$
is measurable with respect to this $\sigma$-algebra.
We define $\mathcal{B}$ to be the smallest $\sigma$-algebra on
$\Lambda$
such that $\pi\dvtx\Lambda\to\bar\Lambda$ and $f\dvtx\Lambda\to
\Lambda$ are measurable.

For any potential function in $F_\theta(\bar\Lambda)$ we obtain a
unique equilibrium state $\bar\mu$. This is an ergodic $\bar{f}$-invariant
probability measure defined on $(\bar\Lambda,\overline{\mathcal{B}})$.
Define $\mu$ on $(\Lambda,\mathcal{B})$ to be the unique
$f$-invariant measure such that $\pi_*\mu=\bar\mu$. Again,
$\mu$ is an ergodic probability measure.

We assume that
there is an integer $m\ge1$ such that all entries of $A^m$ are nonzero.
Then the shift $f$ is mixing with respect to $\mu$.

\begin{pf*}{Proof of Theorem~\ref{thmm-WW2}}
To each $y\in\bar\Lambda$ associate a $y^*\in\Lambda$ such that
(i)~$y^*_i=y_i$ for all $i\ge0$ and
(ii) $x_0=y_0$ implies that $x^*_i=y^*_i$ for each $i\le0$
(e.g., for the full shift, take
$y^*_i=0$ for $i<0$).

Given the observable $v\in F_\theta(\Lambda)$,
define $\chi_1(x)=\sum_{n=0}^\infty v(f^n x^*)-v(f^nx)$.
Then $\chi_1\in L^\infty$ and $v=\hat v+\chi_1\circ f-\chi_1$
where $\hat v$ ``depends only on the future''
and projects down to an observable $\bar v\dvtx\bar\Lambda\to\R$.
Moreover, by Sinai \cite{Sinai72}, $\bar v \in F_{\theta^{1/2}}(\bar
\Lambda)$.
It is standard that there exist constants $a,C>0$ such that $|\int_{\Lambda}\bar v w\circ f^n \,d\mu|\le C\|\bar v\|_{\theta^{1/2}}\|w\|
_1 e^{-an}$ for all $w\in L^1$, $n\ge1$.
By Proposition~\ref{prop-decomp}, \eqref{eq-decomp} holds for all $p$
(even $p=\infty$).
That is, there exist $\bar m,\bar\chi_2\in L^\infty(\bar\Lambda)$
such that
$\bar v=\bar m+\bar\chi_2\circ\bar f-\bar\chi_2$ where
$\bar m\in\ker P$.
It follows that $\hat v=m+\chi_2$ where $m=\bar m\circ\pi$, $\chi
_2=\bar
\chi_2\circ\pi$. Setting $\chi=\chi_1+\chi_2$, we obtain
an $L^\infty$ martingale-coboundary decomposition for $v$ in the sense
of \eqref{eq-decomp_inv}.
Now apply Theorem~\ref{thmm-WWmap_inv}.
\end{pf*}

Our results hold for also for the class of nonuniformly hyperbolic
diffeomorphisms studied by Young \cite{Young98}.
The maps in \cite{Young98} enjoy exponential decay of correlations for
H\"older observables.

More generally, it is possible to consider the situation of Young \cite
{Young99} where the decay of correlations is at a polynomial rate
$n^{-\tau}$.
Provided $\tau>2$ and there is exponential contraction along stable manifolds,
then the conclusion of Theorem~\ref{thmm-WW2} goes through unchanged.
These conditions can be relaxed further; see Section~\ref{sec-gen}.

\section{Iterated WIP for flows}
\label{sec-flow}

In this section, we prove a continuous time version of the iterated WIP
by reducing from continuous time to discrete time.
Theorem~\ref{thmm-WWflow} below is formulated in a purely probabilistic
setting, extending the approach in \cite{MT04,Gouezel07,MZapp}.

We suppose that $f\dvtx\Lambda\to\Lambda$ is a map with ergodic invariant
probability measure $\mu$.
Let $r\dvtx\Lambda\to\R^+$ be an integrable roof function with
$\bar r=\int_\Lambda r \,d\mu$. We suppose throughout that $r$ is bounded below
(away from zero).
Define the suspension
$\Lambda^r=\{(x,u)\in\Lambda\times\R\dvtx0\le u\le r(x)\}/\sim$ where
$(x,r(x))\sim(fx,0)$. Define the suspension flow $\phi
_t(x,u)=(x,u+t)$ computed
modulo identifications. The measure $\mu^r=\mu\times{\rm
Lebesgue}/\bar
r$ is an ergodic invariant probability measure for $\phi_t$.
Using the notation of the \hyperref[sec-intro]{Introduction}, we write
$(\Omega,\nu)=(\Lambda^r,\mu^r)$.

Now suppose that $v\dvtx\Omega\to\R^e$ is integrable with
$\int_{\Omega}v \,d\nu=0$.
Define
the smooth processes $W_n\in C([0,\infty),\R^e)$, $\BBW_n\in
C([0,\infty
),\R^e\times\R^{e\times e})$,
\begin{eqnarray*}
W_n(t) & =&n^{-{1}/2}\int_0^{nt}v
\circ\phi_s \,ds,
\\
\BBW_n^{\beta\gamma}(t) & =&\int_0^t
W_n^\beta \,dW_n^\gamma
=n^{-1}\int_0^{nt}\int
_0^s v^\beta\circ\phi_r
v^\gamma\circ\phi_s \,dr \,ds.
\end{eqnarray*}

Define $\tildev\dvtx\Lambda\to\R^e$ by setting $\tildev(x)=\int_0^{r(x)}v(x,u) \,du$,
and define the cadlag processes
$\tildeW_n\in D([0,\infty),\R^e)$, $\tildeBBW_n\in D([0,\infty),\R
^{e\times e})$,
\begin{eqnarray*}
\tildeW_n(t)&=&n^{-{1}/2}\sum_{j=0}^{[nt]-1}
\tildev\circ f^j,\\
 \tildeBBW_n^{\beta\gamma}(t)&=&\int
_0^t\tildeW_n^\beta d
\tildeW _n^\gamma =n^{-1} \sum
_{0\le i<j\le[nt]-1} \tildev ^\beta\circ f^i
\tildev^\gamma\circ f^j.
\end{eqnarray*}
We assume that the discrete time case is understood, so we have that
%
\begin{equation}
\label{eq-discrete} (\tildeW_n,\tildeBBW_n)
\to_w (\tildeW,\tildeBBW)\qquad\mbox{in } D\bigl([0,\infty),
\R^e\times\R^{e\times e}\bigr),
\end{equation}
where $\tildeW$ is $e$-dimensional Brownian motion and
$\tildeBBW^{\beta\gamma}(t)=\int_0^t\tildeW^\beta \,d\tildeW^\gamma+
\tilde E^{\beta\gamma}t$.
Here, the probability space for the processes on the left-hand side is
$(\Lambda,\mu)$.

Define $H\dvtx\Omega\to\R^e$ by setting $H(x,u)=\int_0^u v(x,s) \,ds$.

\begin{thmm} \label{thmm-WWflow}
Suppose that $\tildev\in L^2(\Lambda)$ and $|H||v|\in L^1(\Omega)$.
Assume \eqref{eq-discrete} and that
%
\begin{eqnarray}
\label{eq-extra}  n^{-1/2}\sup_{t\in[0,T]}|H\circ
\phi_{nt}|&\to_w&0\qquad \mbox{in $C\bigl([0,\infty),
\R^e\bigr)$},
\\
\label{eq-PP}  \lim_{n\to\infty}n^{-1} \biggl\|\max
_{1\le k\le nT} \biggl|\sum_{1\le
i\le k}\tildev\circ
f^i \biggr| \biggr\|_2&=&0.
\end{eqnarray}
Then $(W_n,\BBW_n)\to_w (W,\BBW)$ in $C([0,\infty),\R^e\times\R
^{e\times
e})$ where
the probability space on the left-hand side is $(\Omega,\nu)$, and
\begin{eqnarray*}
 W&=&(\bar r)^{-1/2}\tildeW,\qquad \BBW^{\beta\gamma}(t)=\int
_0^t W^\beta \,dW^\gamma+E^{\beta\gamma}t,
\\
 E^{\beta\gamma}&=&(\bar r)^{-1}\tilde E^{\beta\gamma} +\int
_{\Omega} H^\beta v^\gamma \,d\nu.
\end{eqnarray*}
\end{thmm}

\begin{rmk} \label{rmk-reg}
The regularity conditions on $\tildev$ and $|H||v|$ are satisfied if
$v\in L^\infty(\Omega,\R^e)$ and
$r\in L^2(\Lambda,\R)$, or if $v\in L^2(\Omega,\R^e)$ and $r\in
L^\infty
(\Lambda,\R)$.
Moreover, assumption \eqref{eq-extra} is satisfied under these conditions
by Proposition~\ref{prop-extra}(b).

If $\tildev$ admits an $L^2$ martingale-coboundary decomposition
\eqref
{eq-decomp_inv}, then condition~\eqref{eq-PP} holds by Burkholder's
inequality \cite{Burkholder73}.
\end{rmk}

In the remainder of this section, we prove Theorem~\ref{thmm-WWflow}.
Recall the notation $v_t=\int_0^t v\circ\phi_s \,ds$,
$\tildev_n=\sum_{j=0}^{n-1} \tildev\circ f^j$,
$r_n=\sum_{j=0}^{n-1} r\circ f^j$.
For $(x,u)\in\Omega$ and $t>0$, we define the lap number
$N(t)=N(x,u,t)\in\N$:
\[
N(t)=\max\bigl\{n\ge0\dvtx r_n(x)\le u+t\bigr\}.
\]
Define $g_n(t)=N(nt)/n$.

\begin{lemma} \label{lem-flow1}
$(\tildeW_n,\tildeBBW_n)\circ g_n\to_w
 ((\bar r)^{-1/2}\tildeW,(\bar r)^{-1}\tildeBBW )$
in $D(([0,\infty),\R^e\times\R^{e\times e})$.
\end{lemma}

\begin{pf}
By \eqref{eq-discrete},\vspace*{1pt} $(\tildeW_n,\tildeBBW_n)\to_w(\tildeW
,\tildeBBW)$
on $(\Lambda,\mu)$.
Extend $(\tildeW_n,\tildeBBW_n)$ to $\Omega$ by setting
$\tildeW_n(x,u)=\tildeW_n(x)$,
$\tildeBBW_n(x,u)=\tildeBBW_n(x)$.

We claim that
$(\tildeW_n,\tildeBBW_n)\to_w (\tildeW,\tildeBBW)$ on $(\Omega
,\nu)$.
Define $\bar g(t)=t/\bar r$.
By the ergodic theorem,
$g_n(t)=N(nt)/n=tN(nt)/(nt)\to\bar g(t)$ almost everywhere on $(\Omega
,\nu)$.
Hence, $(\tildeW_n,\tildeBBW_n,g_n)\to_w(\tildeW,\tildeBBW,\bar
g)$ on
$(\Omega,\nu)$.
It follows from the continuous mapping theorem that
\begin{eqnarray*}
\bigl\{(\tildeW_n,\tildeBBW_n)\circ g_n(t),
t\ge0\bigr\} & \to_w& \bigl\{(\tildeW,\tildeBBW)\circ g(t), t\ge0
\bigr\}
\\
& =&\bigl\{ \bigl(\tildeW(t/\bar r),\tildeBBW(t/\bar r)\bigr), t\ge0\bigr\}
\\
& =&\bigl\{ \bigl((\bar r)^{-1/2}\tildeW(t),(\bar r)^{-1}
\tildeBBW(t)\bigr), t\ge0\bigr\}
\end{eqnarray*}
on $(\Omega,\nu)$ completing the proof.

It remains to verify the claim,
Let $c=\operatorname{ess inf} r$
and form the probability space $(\Omega,\mu_c)$
where $\mu_c=(\mu\times{\rm Lebesgue}|_{[0,c]})/c$.
Then it is immediate that $(\tildeW_n,\tildeBBW_n)\to_w (\tildeW
,\tildeBBW)$ on $(\Omega,\mu_c)$.
To pass from $\mu_c$ to $\nu$, and hence to prove the claim, we
apply \cite{Zweimueller07}, Theorem~1.
Since $\mu_c$ is absolutely continuous with respect to $\nu$, it
suffices to prove for all
$\varepsilon,T>0$ that
%
\begin{equation}
\label{eq-P} \lim_{n\to\infty}\mu_r \Bigl(\sup
_{t\in[0,T]}\bigl|P_n(t)\circ f-P_n(t)\bigr|>
\varepsilon \Bigr)=0,
\end{equation}
for $P_n=\tildeW_n$ and $P_n=\tildeBBW_n$.
We give the details for the latter since that is the more complicated case.
Compute that
$\tildeBBW_n^{\beta\gamma}(t)\circ f-\tildeBBW_n^{\beta\gamma
}(t)=n^{-1} \sum_{1\le i<[nt]}\tildev^\gamma\circ f^i \tildev^\beta
\circ f^{[nt]}
-n^{-1} \sum_{1\le j<[nt]}\tildev^\gamma \tildev^\beta\circ f^j$
and so
\begin{eqnarray*}
\Bigl\|\sup_{[0,T]}\bigl|\tildeBBW_n^{\beta\gamma}\circ f-
\tildeBBW _n^{\beta
\gamma}\bigr| \Bigr\|_1 & \le&\bigl\|
\tildev^\beta\bigr\|_2 n^{-1} \biggl\|\max
_{1\le
k\le nT} \biggl|\sum_{1\le i<k}
\tildev^\gamma\circ f^i \biggr|\biggr \|_2
\\
& &{}+ \bigl\|\tildev^\gamma\bigr\|_2 n^{-1}\biggl \|\max
_{1\le k\le
nT} \biggl|\sum_{1\le j\le k}
\tildev^\beta\circ f^j \biggr| \biggr\|_2 \to0
\end{eqnarray*}
by \eqref{eq-PP}.
Hence, \eqref{eq-P} follows from Markov's inequality.
\end{pf}

It follows from the definition of lap number that
\[
\phi_t(x,u)=\bigl(f^{N(t)}x,u+t-r_{N(t)}(x)\bigr).
\]
We have the decomposition
%
\begin{eqnarray}
\label{eq-vV}
v_t(x,u) & =&\int_0^{N(t)}v
\bigl(\phi_s(x,0)\bigr) \,ds+H\circ\phi_t(x,u)-H(x,u)
\nonumber
\\[-8pt]
\\[-8pt]
\nonumber
& =&\tildev_{N(t)}(x)+H\circ\phi_t(x,u)-H(x,u).
\end{eqnarray}

We also require the following elementary result.

\begin{prop}\label{prop-elem} Let $a_n$ be a real sequence and $b>0$.
If $\lim_{n\to\infty}n^{-b}a_n=0$, then $\lim_{n\to\infty
}n^{-b}\sup_{t\in[0,T]}|a_{[nt]}|=0$.
\end{prop}

\begin{lemma} \label{lem-flow2}
$(W_n,\BBW_n) =(\tildeW_n,\tildeBBW_n)\circ g_n + F_n$,
where
$F_n\to_w F$ in $D([0,\infty),\break \R^e\times\R^{e\times e})$ and
$F(t)=
(0,\int_\Omega H^\beta v^\gamma \,d\nu )t$.
\end{lemma}

\begin{pf}
Using \eqref{eq-vV}, we can write
\[
W_n(t)=n^{-1/2}v_{nt}=n^{-1/2}
\tildev_{N(t)}+n^{-1/2}H\circ\phi _{nt}-n^{-1/2}H.
\]
By definition, $\tildeW_n(N(nt)/n)=n^{-1/2}\tildev_{N(t)}$.
Hence, by assumption \eqref{eq-extra}, we obtain
the required decomposition for $W_n$.

Similarly,
%
\begin{eqnarray}
\label{eq-wV}
\BBW^{\beta\gamma}_n(t) & =&\int
_0^t W_n^\beta
\,dW_n^\gamma = \int_0^t
v^\beta_{ns} v^\gamma\circ\phi_{ns} \,ds
\nonumber
\\[-8pt]
\\[-8pt]
\nonumber
& = &\int_0^t \bigl[\tildev^\beta_{N(ns)}+H^\beta
\circ\phi _{ns}-H^\beta\bigr] v^\gamma\circ
\phi_{ns} \,ds = A_n(t)+B_n(t),
\end{eqnarray}
where
\[
A_n(t)=\int_0^t
\tildev^\beta_{N(ns)}v^\gamma\circ\phi_{ns}
\,ds,\qquad B_n(t)=n^{-1}\int_0^{nt}
\bigl[H^\beta\circ\phi_s-H^\beta\bigr]
v^\gamma \circ\phi _s \,ds.
\]

By the ergodic theorem,
\[
n^{-1}H^\beta\int_0^n
v^\gamma\circ\phi_s \,ds =H^\beta(n)^{-1}
\int_0^nv^\gamma\circ
\phi_s \,ds \to H^\beta\int_{\Omega}v^\gamma
\,d\nu=0.
\]
Hence, by Proposition~\ref{prop-elem},
$n^{-1}\sup_{t\in[0,T]}|H^\beta\int_0^{nt} v^\gamma\circ\phi_s
\,ds|\to
0$ a.e.
Similarly,
\[
n^{-1}\int_0^n H^\beta
\circ\phi_s v^\gamma\circ\phi_s \,ds
=n^{-1}\int_0^n\bigl(H^\beta
v^\gamma\bigr) \circ\phi_s \,ds \to\int_{\Omega}H^\beta
v^\gamma \,d\nu.
\]
Applying Proposition~\ref{prop-elem} with $b=1$ and $a_n=\int_0^n
H^\beta\circ\phi_s v^\gamma\circ\phi_s \,ds-\break n\int_{\Omega}H^\beta
v^\gamma \,d\nu$, we obtain that
$n^{-1}\int_0^{nt} H^\beta\circ\phi_s v^\gamma\circ\phi_s \,ds\to
\int_{\Omega}H^\beta v^\gamma \,d\nu$
uniformly on $[0,T]$ a.e.
Hence, $B_n(t)\to t\int_{\Omega}H^\beta v^\gamma \,d\nu$
uniformly on $[0,T]$ a.e.

To deal with the term $A_n$, we
introduce the return times $t_{n,j}=t_{n,j}(x,u)$, with
$0=t_{n,0}<t_{n,1}<t_{n,2}<\cdots$ such that $N(nt)=j$ for $t\in
[t_{n,j},t_{n,j+1})$.
Note that $t_{n,j}(x,u)=(r_j(x)-u)/n$ for $j\ge1$.
Since $r$ is bounded below, we have that $\lim_{j\to\infty
}t_{n,j}=\infty$
for each $n$.

Compute that
\begin{eqnarray*}
A_n(t) & =&\sum_{j=0}^{N(nt)-1}\int
_{t_{n,j}}^{t_{n,j+1}}\tildev_j^\beta
v^\gamma\circ\phi_{ns} \,ds +\int_{t_{n,N(nt)}}^t
\tildev_{N(nt)}^\beta v^\gamma\circ\phi_{ns}
\,ds
\\
& =&\sum_{j=0}^{N(nt)-1}\tildev_j^\beta
\int_{t_{n,j}}^{t_{n,j+1}}v^\gamma\circ
\phi_{ns} \,ds +\tildev_{N(nt)}^\beta\int
_{t_{n,N(nt)}}^t v^\gamma\circ
\phi_{ns} \,ds.
\end{eqnarray*}
For $j\ge1$,
\begin{eqnarray*}
\int_{t_{n,j}}^{t_{n,j+1}}v\circ\phi_{ns} \,ds & =&
\int_{t_{n,j}}^{t_{n,j+1}}v\bigl(f^jx,u+ns-r_j(x)
\bigr) \,ds
\\
& =&n^{-1}\int_0^{r(f^jx)}v
\bigl(f^jx,s\bigr) \,ds =n^{-1}\tildev\circ f^j,
\end{eqnarray*}
and similarly we can write
$\int_0^{t_{n,1}}v\circ\phi_{ns} \,ds=n^{-1}\int_u^{r(x)}v(x,s)
\,ds=n^{-1}\tildev+O(1/n)$ a.e.

By definition, $\tildeBBW_n(N(nt)/n)=n^{-1}\sum_{j=0}^{N(nt)-1}\tildev
_j \tildev\circ f^j$.
Hence, we have shown that
$A_n(t)=\tildeBBW_n\circ g_n(t)+C_n(t)+O(1/n)$ a.e., where
$C_n^{\beta\gamma}(t)=\tildev_{N(nt)}^\beta\int_{t_{n,N(nt)}}^t
v^\gamma
\circ\phi_{ns} \,ds$.

Finally, we note that
\begin{eqnarray*}
 \int_{t_{n,N(nt)}}^t v\circ\phi_{ns} \,ds &=&
\int_{t_{n,N(nt)}}^t v\bigl(f^{N(nt)}x,u+ns-r_{N(nt)}(x)
\bigr) \,ds
\\
& = &n^{-1}\int_0^{u+t-r_{N(nt)}(x)}v
\bigl(f^{N(nt)}x,s\bigr) \,ds\\
&=&n^{-1}H\bigl(f^{N(nt)}x,u+t-r_{N(nt)}(x)
\bigr)
\\
& =&n^{-1}H\circ\phi_{nt}.
\end{eqnarray*}
Hence, $C_n^{\beta\gamma}=\tildeW_n^\beta\circ g_n(t)\cdot
n^{-1/2}H^\gamma\circ\phi_{nt}\to_w0$ by
Lemma~\ref{lem-flow1} and assumption~\eqref{eq-extra}.
\end{pf}

\begin{pf*}{Proof of Theorem~\ref{thmm-WWflow}}
This is immediate from Lemmas \ref{lem-flow1} and \ref{lem-flow2}.
\end{pf*}

\begin{prop} \label{prop-extra}
Sufficient conditions for assumption \eqref{eq-extra}
to hold are that
\textup{(a)} $H\in L^{2+}(\Omega,\R^e)$, or
\textup{(b)} $\tildev_*\in L^2(\Lambda)$,
where $\tildev_*(x)=\int_0^{r(x)}|v(x,u)| \,du$.
\end{prop}

\begin{pf}
In both cases, we prove that
$n^{-1/2}H\circ\phi_n\to0$ a.e. By Proposition~\ref{prop-elem},
$\sup_{t\in[0,T]}H\circ\phi_{nt}\to0$ a.e.

(a) Choose $\delta>0$ such that $H\in L^{2+\delta}$ and $\tau
<\frac{1}2$ such that
\mbox{$\tau(2+\delta)>1$}.
Since $\|H\circ\phi_n\|_{2+\delta}=\|H\|_{2+\delta}$, it follows from
Markov's inequality that $\nu(|H\circ\phi_n|>n^\tau)\le\|H\|
_{2+\delta
}n^{-\tau(2+\delta)}$ which is summable.
By Borel--Cantelli, there is a constant $C>0$ such that
$|H\circ\phi_n|\le Cn^{-\tau}$ a.e., and hence $n^{-1/2}H\circ\phi
_n\to
0$ a.e.

(b)
Since $\tildev_*^2\in L^1(\Lambda)$, it follows from the ergodic
theorem that
$n^{-1/2}\tildev_*\circ f^n\to0$ a.e.
Moreover, $N(nt)/n\to1/\bar r$ a.e.
on $(\Omega,\nu)$ and hence
$n^{-1/2}\tildev_*\circ f^{[N(nt)]}\to0$ a.e.
The result follows since $|H(x,u)|\le\tildev_*(x)$ for all $x,u$.
\end{pf}

\begin{rmk} \label{rmk-extra}
The sufficient conditions in Proposition~\ref{prop-extra} imply almost
sure convergence, uniformly on $[0,T]$, for the term $F_n$ in Lemma~\ref
{lem-flow2}.
\end{rmk}

\section{Moment estimates}
\label{sec-moment}

In this section, we obtain some moment estimates that are required to
apply rough path theory.
(Proposition~\ref{prop-moment_flow} below is also required for part of
Theorem~\ref{thmm-WW}; see the proof of Corollary~\ref{cor-WWflowmix}.)

\subsection{Discrete time moment estimates}

Let $f\dvtx\Lambda\to\Lambda$ be a map (invertible or
noninvertible) with
invariant probability measure $\mu$.
Suppose that
$v\dvtx\Lambda\to\R^e$ is a mean zero observable lying in $L^\infty$.
Define
\[
v_n=\sum_{j=0}^{n-1}v\circ
f^j,\qquad S_n^{\beta\gamma}= \sum
_{0\le i<j<n}v^\beta\circ f^i
v^\gamma\circ f^j.
\]

\begin{prop} \label{prop-moment}
Suppose that $v\dvtx\Lambda\to\R^e$ lies in $L^\infty$ and admits
an $L^p$
martingale-coboundary decomposition \eqref{eq-decomp_inv} for some
$p\ge3$.
Then there exists a constant $C>0$ such that
\[
\Bigl\|\max_{0\le j\le n}|v_j| \Bigr\|_{2p}\le
Cn^{1/2}, \qquad\Bigl\| \max_{0\le j\le n}|S_j|
\Bigr\|_{2p/3} \le Cn \qquad\mbox{for all $n\ge1$.}
\]
\end{prop}

\begin{pf}
The estimate $\|v_n\|_{2p}\ll n^{1/2}$ is proved in \cite{MN08},
equation (3.1).
Since $v_{n+a}-v_a=_d v_n$ for all $a,n$, the result for
$\max_{0\le j\le n}|v_j|$ follows by \cite{Serfling70}, Corollary B1
(cf. \cite{MTorok12}, Lemma~4.1).

To estimate $S_n$ write
\[
S_n^{\beta\gamma}= \sum_{0\le i<j<n}m^\beta
\circ f^i v^\gamma\circ f^j+ \sum
_{1\le j<n}\bigl(\chi^\beta\circ f^j-
\chi^\beta\bigr) v^\gamma\circ f^j.
\]
We have
$\|\sum_{1\le j<n}\chi^\beta\circ f^j v^\gamma\circ f^j\|_p\le
n\|\chi^\beta v^\gamma\|_p\le n\|\chi^\beta\|_p\|v^\gamma\|
_\infty$ and\break
$\|\sum_{1\le j<n}\chi^\beta v^\gamma\circ f^j\|_p\le
\|\chi^\beta\|_p
\|\sum_{1\le j<n}v^\gamma\circ f^j\|_\infty\le
n\|\chi^\beta\|_p\|v^\gamma\|_\infty$.

Next, we estimate
$I_n=\sum_{0\le i<j<n}m^\beta\circ f^i v^\gamma\circ f^j$.
Passing to the natural extension $\tilde f\dvtx\tilde\Lambda\to
\tilde\Lambda$
in the noninvertible case (and taking $\tilde f=f$ in the invertible case),
we have
\[
\tilde I_n=\sum_{0\le i<j<n}\tilde
m^\beta\circ\tilde f^i \tilde v^\gamma\circ\tilde
f^j= \biggl(\sum_{-n\le i<j<0} \tilde
m^\beta\circ\tilde f^i \tilde v^\gamma\circ\tilde
f^j \biggr)\circ\tilde f^n= \tilde I_n^-\circ
\tilde f^n,
\]
so we reduce to estimating
$\tilde I_n^-=\sum_{-n\le i<j<0} \tilde v^\gamma\circ\tilde f^j
\tilde
m^\beta\circ\tilde f^i$.

Now,
\[
\tilde I_n^-=\sum_{k=1}^n
X_k \qquad\mbox{where } X_k= \biggl(\sum
_{-k<j<0}\tilde v^\gamma\circ\tilde f^j
\biggr)\tilde m^\beta\circ\tilde f^{-k}.
\]
Recall that $E(\tilde m^\beta\circ\tilde f^i|\tilde f^{-i-1}\tilde
{\mathcal{B}})=0$.
Hence
$E(X_k|\tilde f^{k-1}\tilde{\mathcal{B}})=0$, and so
$\{X_k; k\ge1\}$ is a sequence of martingale differences.
For $p'>1$, Burkholder's inequality \cite{Burkholder73} states that
$\|\tilde I_n^-\|_{p'}
\ll\|(\sum_{k=1}^nX_k^2)^{1/2}\|_{p'}$,
and it follows for $p'\ge2$ that
%
\begin{equation}
\label{eq-Burk} \bigl\|\tilde I_n^-\bigr\|_{p'}^2 \ll
\sum_{k=1}^n\|X_k
\|_{p'}^2.
\end{equation}
Taking $p'=2p/3$, it follows from H\"older's inequality that
\[
\|X_k\|_{2p/3}\le \biggl\|\sum_{-k<j<0}
\tilde v^\gamma\circ\tilde f^j\biggr \|_{2p}\bigl\| \tilde
m^\beta\circ\tilde f^{-k}\bigr\|_p =
\bigl\|v^\gamma_{k-1}\bigr\|_{2p} \bigl\|m^\beta\bigr\|
_p\ll k^{1/2}.
\]
Hence, $\|\tilde I_n^-\|_{2p/3}\ll n$ and so
$\|S_n\|_{2p/3}\ll n$.

This time we cannot apply the maximal inequality of \cite{Serfling70}
since we do not have a good estimate for $S_{a+n}-S_a$ uniform in $a$.
However, we claim that $\|S_{a+n}-S_a\|_{2p/3}\ll
n+n^{1/2}a^{1/2}$. Set
$A_{a,n}=(\sum_{k=a+1}^{a+n}b_k^2)^{1/2}$ with $b_k=k^{1/2}$.
By the claim,
$\|S_{a+n}-S_a\|_{2p/3}\ll
A_{a,n}$ and it follows from \cite{Moricz76}, Theorem A (see also
references therein) that
$\|\max_{0\le j\le n}|S_j| \|_{2p/3}\ll n$ as required.

For the claim, observe that
\begin{eqnarray*}
S^{\beta\gamma}_{a+n}-S^{\beta\gamma}_a & =&\sum
_{j=a}^{a+n-1}\sum_{i=0}^{j-1}v^\beta
\circ f^i v^\gamma\circ f^j
\\
& =&\sum_{j=a}^{a+n-1}\sum
_{i=0}^{a-1}v^\beta\circ f^i
v^\gamma \circ f^j +\sum_{j=a}^{a+n-1}
\sum_{i=a}^{j-1}v^\beta\circ
f^i v^\gamma\circ f^j
\\
& =&v_a^\beta v_n^\gamma\circ
f^a+S_n^{\beta\gamma}\circ f^a.
\end{eqnarray*}
Hence,
\[
\bigl\|S^{\beta\gamma}_{a+n}-S^{\beta\gamma}_a
\bigr\|_q\le\bigl\|v_n^\gamma\bigr\| _{2q}\bigl\|
v_a^\beta\bigr\|_{2q}+\|S_n
\|_q \ll n^{1/2}a^{1/2}+n,
\]
for $q=2p/3$. This proves the claim.
\end{pf}

\begin{rmk}
The proof of Proposition~\ref{prop-moment} makes essential use of the
fact that $v\in L^\infty$ \cite{LesigneVolny01,MN08,MTorok12}. Under
this assumption,
the estimate for $\max_{0\le j\le n}|v_j|$ requires only that $p\ge1$
and is optimal in the sense that there are examples where $\lim_{n\to
\infty}\|n^{-1/2}v_n\|_q=\infty$ for all $q>2p$; see \cite{MTorok12},
Remark~3.7.

We conjecture that the optimal estimate for $\max_{0\le j\le n}|S_j|$
is that\break
$\|\max_{0\le j\le n}|S_j| \|_p\ll n$ (for $p\ge2$).
Then we would only require $p>3$ instead of $p>9/2$ in our main results.
\end{rmk}

Recall that
$W_n(t)=n^{-1/2}\sum_{j=0}^{[nt]-1}v\circ f^j$ and
$\BBW_n^{\beta\gamma}(t)=\int_0^t W_n^\beta \,dW_n^\gamma$.
We define the increments
\[
W_n(s,t)=W_n(t)-W_n(s) \quad\mbox{and}\quad
\BBW^{\beta\gamma
}_n(s,t)=\int_s^t
W^\beta_n(s,r) \,dW^\gamma_n(r).
\]

\begin{cor} \label{cor-moment}
Suppose that $v\dvtx\Lambda\to\R^e$ lies in $L^\infty$ and admits
an $L^p$
martingale-coboundary decomposition \eqref{eq-decomp_inv} for some
$p\ge3$.
Then there exists a constant $C>0$ such that
\begin{eqnarray*}
\bigl\|W_n(j/n,k/n)\bigr\|_{2p} &\leq& C\bigl(|k-j|/n\bigr)^{1/2}\quad
\mbox{and}\\
 \bigl\|\BBW _n(j/n,k/n) \bigr\|_{2p/3} &\leq& C|k-j|/n,
\end{eqnarray*}
for all $j,k,n\ge1$.
\end{cor}

\begin{pf}
Let $t>s>0$. By definition,
\begin{eqnarray*}
W_n(s,t) & =&n^{-1/2}\sum_{i=[ns]}^{[nt]-1}v
\circ f^i =n^{-1/2} \Biggl(\sum_{i=0}^{[nt]-[ns]-1}v
\circ f^i \Biggr)\circ f^{[ns]}
\\
& =_d& n^{-1/2}\sum_{i=0}^{[nt]-[ns]-1}v
\circ f^i= n^{-1/2}v_{[nt]-[ns]}.
\end{eqnarray*}
By Proposition~\ref{prop-moment}, assuming without loss that $j<k$,
\[
\bigl\|W_n(j/n,k/n)\bigr\|_{2p}=n^{-1/2}\|v_{k-j}
\|_{2p}\le C\bigl((k-j)/n\bigr)^{1/2}.
\]

Similarly,
\begin{eqnarray*}
\BBW_n(s,t) & =& n^{-1} \sum_{[ns]\le i<j\le[nt]-1}
v^\beta\circ f^i v^\gamma\circ f^j\\
&=&n^{-1} \biggl(\sum_{0\le i<j<[nt]-[ns]-1}
v^\beta \circ f^i v^\gamma\circ f^j
\biggr)\circ f^{[ns]}
\\
& =_d&n^{-1} \sum_{0\le i<j<[nt]-[ns]-1}
v^\beta\circ f^i v^\gamma\circ
f^j=n^{-1}S^{\beta\gamma
}_{[nt]-[ns]}.
\end{eqnarray*}
By Proposition~\ref{prop-moment},
\[
\bigl\|\BBW_n(j/n,k/n)\bigr\|_{2p/3}= n^{-1}
\|S_{k-j}\|_{2p/3}\le C(k-j)/n,
\]
as required.
\end{pf}

\subsection{Continuous time moment estimates}

Let $\phi_t\dvtx\Omega\to\Omega$ be a suspension flow as in
Section~\ref{sec-flow},
with Poincar\'e map $f\dvtx\Lambda\to\Lambda$. As before, we write
$\Omega=\Lambda^r$, $\nu=\mu^r$, where $r\dvtx\Lambda\to\R$ is
a roof
function with
$\bar r=\int r \,d\mu$.
Let $v\dvtx\Omega\to\R^e$ with $\int_\Omega v \,d\nu=0$.

As before, we suppose that $r$ is bounded away from zero, but now
we suppose in addition that $v$ and $r$ lie in $L^\infty$.
(These assumptions can be relaxed, but then the assumption on $p$ has
to be strengthened in the subsequent results.)

Define
\[
v_t=\int_0^t v\circ
\phi_s \,ds, \qquad S_t^{\beta\gamma}= \int_0^t
\int_0^u v^\beta\circ
\phi_s v^\gamma \circ \phi_u \,ds \,du.
\]
Let $\tildev\dvtx\Lambda\to\R^e$ be given by $\tildev(x)=\int_0^{r(x)}v(x,u) \,du$
(so $\tildev$ coincides with the function defined in Section~\ref{sec-flow}).
The assumptions on $v$ and $r$ imply that $\tildev\in L^\infty
(\Lambda
,\mu)$.

\begin{prop} \label{prop-lap}
$N(t)\le[t/\essinf r]+1$ for all
$(x,u)\in\Omega$, $t\ge0$.
\end{prop}

\begin{pf}
Compute that
\begin{eqnarray*}
r_{[t/\essinf r]+2}(x)&=&r(x)+ r_{[t/\essinf r]+1}(fx)\\
&\ge& u+\bigl([t/\essinf r]+1
\bigr)\essinf f> u+t.
\end{eqnarray*}
Hence, the result follows from the definition of lap number.
\end{pf}

\begin{prop} \label{prop-moment_flow}
Suppose that $\tildev\dvtx\Lambda\to\R^e$ admits an $L^p$
martingale-coboundary decomposition \eqref{eq-decomp_inv} for some
$p\ge3$.
Then there exists a constant $C>0$ such that
\[
\|v_t\|_{2p}\le Ct^{1/2},\qquad \|S_t
\|_{2p/3} \le Ct,
\]
for all $t\ge0$.
\end{prop}

\begin{pf}
If $t\le1$, then we have the almost sure estimates
$|v_t|\le\|v\|_\infty t\le\|v\|_\infty t^{1/2}$ and
$|S_t|\le\|v\|_\infty^2t^2\le\|v\|_\infty^2t$. Hence, in the
remainder of the proof, we can suppose that $t\ge1$.

For the $v_t$ estimate, we follow the argument used in \cite{MTorok12},
Lemma~4.1.
By \eqref{eq-vV},
\[
v_t=\tildev_{N(t)}+G(t),
\]
where $G(t)(x,u)=H\circ\phi_t(x,u)-H(x,u)=\int_0^u v(\phi_t(x,s)) \,ds-
\int_0^u v(x,s) \,ds$.
In particular, $\|G(t)\|_\infty\le2\|r\|_\infty\|v\|_\infty
\le2\|r\|_\infty\|v\|_\infty t^{1/2}$.
By Proposition~\ref{prop-lap}, there is a constant $R>0$ such that
$N(t)\le Rt$ for all $t\ge1$.
Hence,
\[
|v_t|\le\max_{0\le j\le Rt}|\tildev_j|+2\|r
\|_\infty\|v\|_\infty t^{1/2}.
\]
By Proposition~\ref{prop-moment},
$\|\max_{0\le j\le Rt}|\tildev_j|\|_{2p}\ll t^{1/2}$.
Since $r$ is bounded above and below, this estimate for
$\max_{0\le j\le Rt}|\tildev_j|$ holds equally in
$L^{2p}(\Lambda)$ and $L^{2p}(\Omega)$.
Hence $\|v_t\|_{2p}\ll t^{1/2}$.

To estimate $S_t$ we make use of decompositions similar to those
in Section~\ref{sec-flow}.
By~\eqref{eq-vV},
\[
S_t^{\beta\gamma} = \int_0^t
v_s^\beta v^\gamma\circ\phi_s \,ds =
\int_0^t \bigl(\tildev_{N(s)}^\beta+G^\beta(s)
\bigr)v^\gamma\circ\phi_s \,ds,
\]
where
$\|\int_0^t G^\beta(s) v^\gamma\circ\phi_s \,ds\|_\infty\le
2|r|_\infty
|v|_\infty^2 t$.
Moreover, in the notation from the proof of Lemma~\ref{lem-flow2} with $n=1$,
\begin{eqnarray*}
\int_0^t \tildev_{N(s)}^\beta
v^\gamma\circ\phi_s \,ds &=& A_1(t)\\
& =& \sum
_{j=0}^{N(t)-1} \tildev^\beta_j
\tildev^\gamma\circ f^j-\tilde v^\beta\int
_0^u v^\gamma\circ\phi_s
\,ds+\tildev^\beta _{N(t)}H^\gamma \circ
\phi_t
\\
& =& \tilde S^{\beta\gamma}_{N(t)}-\tilde v^\beta\int
_0^u v^\gamma \circ
\phi_s \,ds+\tildev^\beta_{N(t)}H^\gamma
\circ\phi_t,
\end{eqnarray*}
where $\tilde S_n$ is as in Proposition~\ref{prop-moment}.
Now $\|\tildev^\beta_{N(t)}\|_\infty\le\|N(t)\|_\infty\|\tildev
^\beta\|
_\infty\le\break  Rt\|r\|_\infty\|v\|_\infty$.
Hence, by Proposition~\ref{prop-moment},
\[
\biggl|\int_0^t \tildev_{N(s)}^\beta
v^\gamma\circ\phi_s \,ds \biggr|\le \max_{j\le Rt}\bigl|
\tilde S_j^{\beta\gamma}\bigr|+(1+Rt)\|r\|_\infty^2
\|v\| _\infty^2 \ll t,
\]
completing the proof.
\end{pf}

Again we recall that
$W_n(t)=n^{-1/2}\int_0^{nt}v\circ\phi_s \,ds$ and
$\BBW_n^{\beta\gamma}(t)=\int_0^t W_n^\beta \,dW_n^\gamma$,
and we define the increments
\[
W_n(s,t)=W_n(t)-W_n(s)\quad\mbox{and}\quad
\BBW^{\beta\gamma
}_n(s,t)=\int_s^t
W^\beta_n(s,r) \,dW^\gamma_n(r).
\]

\begin{cor} \label{cor-moment_flow}
Suppose that $\tildev$ admits an $L^p$ martingale-coboundary
decomposition \eqref{eq-decomp_inv} for some $p\ge3$.
Then there exists a constant $C>0$ such that
\[
\bigl\|W_n(s,t) \bigr\|_{2p} \leq C|t-s|^{1/2} \quad\mbox{and}\quad
\bigl\|\BBW _n(s,t) \bigr\| _{2p/3} \leq C|t-s|,
\]
for all $s,t\ge0$.
\end{cor}

\begin{pf}
This is almost identical to the proof of Corollary~\ref{cor-moment}.
\end{pf}

\begin{rmk}
Any hyperbolic basic set for an Axiom A flow can
be written as a suspension over a mixing hyperbolic basic set $f\dvtx
\Lambda
\to\Lambda$ with a H\"older roof function $r$.
Since every H\"older mean zero observable $\tilde v\dvtx\Lambda\to\R
^e$ admits
an $L^\infty$ martingale-coboundary decomposition, it follows that
Proposition~\ref{prop-moment_flow}
and Corollary~\ref{cor-moment_flow} hold for all $p$.
\end{rmk}

\section{Applications of Theorem \texorpdfstring{\protect\ref{thmm-WWflow}}{6.1}}
\label{sec-flowpf}

In this section, we apply Theorem~\ref{thmm-WWflow} to a large class of
uniformly and nonuniformly hyperbolic flows. In particular, we complete
the proof of
Theorem~\ref{thmm-WW}.
Our main results \,do not require mixing assumptions on the flow, but the
formulas simplify in the mixing case.

Let $\phi_t\dvtx\Omega\to\Omega$ be a suspension flow as in
Section~\ref{sec-flow},
with mixing Poincar\'e map $f\dvtx\Lambda\to\Lambda$. As before, we write
$\Omega=\Lambda^r$, $\nu=\mu^r$, where $r\dvtx\Lambda\to\R$ is
a roof
function with
$\bar r=\int r \,d\mu$.

\textit{Nonmixing flows.}
First, we consider the case where $\phi_t$ is not mixing. (As usual, we
suppose that the Poincar\'e map $f$ is mixing.)

\begin{cor} \label{cor-WWflow}
Suppose that $f\dvtx\Lambda\to\Lambda$ is mixing and that $r\in
L^1(\Lambda
)$ is bounded away from zero.
Let $v\in L^1(\Omega,\R^e)$ with $\int_\Omega v \,d\nu=0$.
Suppose further that $|H||v|$ is integrable and that assumption~\eqref
{eq-extra} is satisfied.

Assume that $\tildev$ admits
a martingale-coboundary decomposition \eqref{eq-decomp_inv} with $p=2$.
Then the conclusion of Theorem~\ref{thmm-WWflow} is valid.
Moreover,
\begin{eqnarray*}
\Sigma^{\beta\gamma}&=&\Cov^{\beta\gamma}W(1)\\
&=&(\bar r)^{-1}\int
_\Lambda \tildev^\beta\tildev^\gamma \,d\mu+(
\bar r)^{-1}\sum_{n=1}^\infty \int
_\Lambda\bigl(\tildev^\beta\tildev^\gamma
\circ f^n+\tildev^\gamma \tildev ^\beta\circ
f^n\bigr) \,d\mu,
\end{eqnarray*}
and $\BBW^{\beta\gamma}(t)=\int_0^t W^\beta\circ \,dW^\gamma+\frac{1}2
D^{\beta\gamma}t$
where
\[
D^{\beta\gamma}=(\bar r)^{-1}\sum_{n=1}^\infty
\int_\Lambda \bigl(\tildev ^\beta\tildev^\gamma
\circ f^n-\tildev^\gamma\tildev^\beta\circ
f^n\bigr) \,d\mu+\int_\Omega\bigl(H^\beta
v^\gamma-H^\gamma v^\beta\bigr) \,d\nu.
\]
\end{cor}

\begin{pf}
By Theorem~\ref{thmm-WWmap}, condition \eqref{eq-discrete} is satisfied.
Specifically,\break $(\tildeW_n,\tildeBBW_n)\to_w(\tildeW,\tildeBBW)$ where
$\tildeW$ is a Brownian motion with
$\Cov^{\beta\gamma}(\tildeW(1))= \int_\Lambda\tildev^\beta
\tildev
^\gamma \,d\mu+\sum_{n=1}^\infty\int_\Lambda(\tildev^\beta\tildev
^\gamma\circ f^n+\tildev^\gamma\tildev^\beta\circ f^n) \,d\mu$ and
$\tildeBBW_n^{\beta\gamma}(t)=\int_0^t \tildeW^\beta \,d\tildeW
^\gamma
+\tilde E^{\beta\gamma}t$.

By Remark~\ref{rmk-reg}, hypothesis \eqref{eq-PP} is satisfied.
Hence, by Theorem~\ref{thmm-WWflow},
$(W_n, \BBW_n)\to_w(W,\BBW)$ where $W=(\bar r)^{-1/2}\tildeW$ and
$\BBW^{\beta\gamma}(t)= \int_0^t W_n^\beta \,dW_n^\gamma+E^{\beta
\gamma}t$.
It is immediate that $\Sigma=\Cov(W(1))$ has the desired form.
Moreover, by Theorems~\ref{thmm-WWmap} and \ref{thmm-WWflow},
\[
E^{\beta\gamma}=(\bar r)^{-1}\sum_{n=1}^\infty
\int_\Lambda\tildev ^\beta \tildev^\gamma
\circ f^n \,d\mu+\int_\Omega H^\beta
v^\gamma \,d\nu.
\]
The Stratonovich correction gives
\begin{eqnarray*}
D^{\beta\gamma} & =&2E^{\beta\gamma}-\Sigma^{\beta\gamma}
\\
& = &(\bar r)^{-1} \Biggl\{\sum_{n=1}^\infty
\int_\Lambda\bigl(\tildev^\beta \tildev^\gamma
\circ f^n -\tildev^\gamma\tildev^\beta\circ
f^n\bigr) \,d\mu-\int_\Lambda\tildev ^\beta
\tildev^\gamma \,d\mu \Biggr\}\\
&&{}+2\int_\Omega
H^\beta v^\gamma \,d\nu.
\end{eqnarray*}

To complete the proof, we show that
$(\bar r)^{-1}\int_\Lambda\tildev^\beta\tildev^\gamma \,d\mu=
\int_{\Omega} H^\beta v^\gamma \,d\nu+
\int_{\Omega} H^\gamma v^\beta \,d\nu$.
Compute that
\begin{eqnarray*}
 \int_\Lambda\tildev^\beta\tildev^\gamma d
\mu &=& \int_\Lambda \biggl\{\int_0^{r(x)}v^\beta(x,u)
\,du \int_0^{r(x)}v^\gamma(x,s) \,ds \biggr
\}\,d\mu
\\
& =& \int_\Lambda \int_0^{r(x)}v^\beta(x,u)
\biggl\{ \int_0^uv^\gamma(x,s) \,ds+
\int_u^{r(x)}v^\gamma(x,s) \,ds \biggr\} \,du
\,d\mu
\\
& =& \int_\Lambda \int_0^{r(x)}v^\beta(x,u)
H^\gamma(x,u) \,du \,d\mu
\\
&&{} + \int_\Lambda \int_0^{r(x)}v^\gamma(x,s)
\biggl(\int_0^s v^\beta(x,u) \,du
\biggr) \,ds \,d\mu
\\
& = &\bar r\int_{\Omega}v^\beta H^\gamma \,d\nu
+ \int_\Lambda \int_0^{r(x)}v^\gamma(x,s)
H^\beta(x,s) \,ds
\\
& = &\bar r\int_{\Omega}v^\beta H^\gamma \,d\nu
+ \bar r\int_{\Omega}v^\gamma H^\beta \,d\nu,
\end{eqnarray*}
as required.
\end{pf}

\begin{rmk} \label{rmk-WWflow}
Corollary~\ref{cor-WWflow} applies directly to H\"older observables of
semiflows that are suspensions of the uniformly and nonuniformly
expanding maps considered in Section~\ref{sec-NUE}, and of flows that
are suspensions of the uniformly and nonuniformly hyperbolic
diffeomorphisms considered in Section~\ref{sec-NUH}. In particular,
this includes Axiom A flows and nonuniformly hyperbolic flows that are
suspensions over Young towers with exponential tails.
\end{rmk}

\textit{Mixing flows.}
Under additional conditions, we obtain the formulas for $\Sigma$ and
$D$ promised in Theorem~\ref{thmm-WW}(b).

\begin{cor} \label{cor-WWflowmix}
Assume the set up of Corollary~\ref{cor-WWflow}.
Suppose further that $v\in L^\infty$, and that $\tildev$ admits
a martingale-coboundary decomposition \eqref{eq-decomp_inv} with $p=3$.
If the integral $\int_0^\infty\int_\Omega v^\beta v^\gamma\circ
\phi_t
\,d\nu \,dt$ exists, then
\[
\Sigma^{\beta\gamma}= \int_0^\infty\int
_\Omega\bigl(v^\beta v^\gamma\circ
\phi_t+ v^\gamma v^\beta\circ\phi_t
\bigr) \,d\nu \,dt
\]
and
\[
D^{\beta\gamma}= \int_0^\infty\int
_\Omega\bigl( v^\beta v^\gamma\circ
\phi_t- v^\gamma v^\beta\circ\phi_t\bigr) \,d
\nu \,dt.
\]
\end{cor}

\begin{pf}
It follows from \cite{Burkholder73} that $\|W_n\|_p=O(1)$, and hence that\break
$\E_\nu|W_n|^q\to\E|W|^q$ for all $q<p$.
In particular, taking $q=2$, we deduce that
$\Cov_\nu(W_n(1))\to\Sigma$.
Moreover, the calculation in the proof of Theorem~\ref{thmm-WWmap}
shows that
$\Sigma^{\beta\gamma}=\int_0^\infty\int_\Omega(v^\beta v^\gamma
\circ
\phi_t+v^\gamma v^\beta\circ\phi_t) \,d\nu \,dt$.
Similarly\break $\E_\nu(\int_0^1 W_n^\beta \,dW_n^\gamma)\to\int_0^\infty
\int_\Omega v^\beta v^\gamma\circ\phi_t \,d\nu \,dt$.

Since $v\in L^\infty$ and $p=3$, it follows from Proposition~\ref
{prop-moment_flow} that
$\|\int_0^1 W_n^\beta \,dW_n^\gamma\|_2=O(1)$.
[In fact, we require only that
$\|\int_0^1 W_n^\beta \,dW_n^\gamma\|_q=O(1)$ for some $q>1$.]
Hence, $\E_\nu(\int_0^1 W_n^\beta \,dW_n^\gamma)\to
E^{\beta\gamma}$, and so
\begin{eqnarray*}
D^{\beta\gamma} & =&2E^{\beta\gamma}-\Sigma^{\beta\gamma}
\\
& =& 2\int_0^\infty\int_\Omega
v^\beta v^\gamma\circ\phi_t \,d\nu \,dt-\int
_0^\infty\int_\Omega
\bigl(v^\beta v^\gamma\circ\phi_t+v^\gamma
v^\beta \circ\phi_t\bigr) \,d\nu \,dt
\\
& =&\int_0^\infty\int_\Omega
\bigl(v^\beta v^\gamma\circ\phi_t-v^\gamma
v^\beta\circ\phi_t\bigr) \,d\nu \,dt,
\end{eqnarray*}
as required.
\end{pf}

\begin{pf*}{Proof of Theorem~\ref{thmm-WW}}
We use the fact that every hyperbolic basic set for an Axiom A flow can
be written as a suspension over a mixing hyperbolic basic set $f\dvtx
\Lambda
\to\Lambda$ with a H\"older roof function $r$.
Any H\"older mean zero observable $v\dvtx\Lambda\to\R^e$ admits
an $L^\infty$ martingale-coboundary decomposition.
Hence Theorem~\ref{thmm-WW} follows from Corollaries \ref{cor-WWflow}
and \ref{cor-WWflowmix}.
\end{pf*}

\section{Smooth approximation theorem}
\label{sec-approx}

In this section, we prove Theorems \ref{thmm-approx} and \ref
{thmm-approx2}. To do so, we need a few tools from \emph{rough path
theory} that
allow us to lift the iterated WIP into a convergence result for
fast-slow systems. We do not need to introduce much new terminology
since the tools we need are to some extent {prepackaged} for our
purposes. For the continuous time results, we use the standard rough
path theory \cite{Lyons98}, but for the discrete time results we use
results of \cite{Kelly14}.

\subsection{Rough path theory in continuous time}

Let $U_n \dvtx[0,T]\to\R^e$ be a path of bounded variation. Then we can
define the iterated integral $\BBU_n \dvtx[0,T]\to\R^{e\times e}$ by
%
\begin{equation}
\label{e:bbu_flow} \BBU_n(t) = \int_0^t
U_n(r) \,dU_n(r),
\end{equation}
where the integral is uniquely defined in the Riemann--Stieltjes sense.
As usual, we define the increments
\[
U_n(s,t) = U_n(t)-U_n(s) \quad\mbox{and}\quad
\BBU_n(s,t) = \int_s^t
U_n(s,r) \,dU_n(r).
\]
Suppose that $a\dvtx\R^d\to\R^d$ is $C^{1+}$ and $b\dvtx\R^d\to
\R^{d\times e}$
is $C^{3}$,
and let $X_n$ be the solution to the equation
%
\begin{equation}
\label{e:cts_RDE} X_n (t) = \xi+ \int_0^t
a\bigl(X_n(s)\bigr) \,ds+ \int_0^t b
\bigl(X_n(s)\bigr) \,dU_n(s),
\end{equation}
which is well-defined for each $n$, and moreover has a unique solution
for every initial condition $\xi\in\R^d$. To characterize the limit of
$X_n$, we use the following standard tool from rough path theory.

\begin{thmm}\label{thmm:cts_sm}
Suppose that $(U_n,\BBU_n) \to_w (U,\BBU)$ in $C([0,\infty),\R
^e\times\R
^{e\times e})$, where $U$ is Brownian motion and where $\BBU$ can be written
%
\begin{equation}
\label{e:strato_E} \BBU(t) = \int_0^t U(s) \circ
dU(s) + D t,
\end{equation}
for some constant matrix $D \in\R^{e\times e}$. Suppose moreover that
there exist $C>0$ and $q>3$ such that
%
\begin{equation}
\label{e:cts_kolm}\bigl \|U_n(s,t) \bigr\|_{2q} \leq
C|t-s|^{1/2} \quad\mbox{and}\quad\bigl \|\BBU _n(s,t) \bigr\|_{q}
\leq C |t-s|,
\end{equation}
hold for all $n\ge1$ and $s,t \in[0,T]$. Then $X_n \to_w X$ in
$C([0,\infty),\R^d)$, where
%
\begin{equation}
\label{e:RDE_SDE} dX = \biggl(a(X) + \sum_{\alpha,\beta,\gamma}D^{\beta\gamma}
\partial^\alpha b^{\beta}(X) b^{ \alpha\gamma}(X) \biggr)\,dt +b(X)
\circ dW.
\end{equation}
If \eqref{e:cts_kolm} holds for all $q<\infty$, then the $C^3$
condition on $b$ can be relaxed to $C^{2+}$.
\end{thmm}

This result has been used in several contexts \cite
{friz09walk,lejay05}, so we only sketch the proof.

\begin{pf*}{Proof of Theorem \ref{thmm:cts_sm}}
First, suppose that $a\in C^{1+}$, $b\in C^{3}$.
By \cite{friz10}, Theorem~12.10, we know that the map $(U_n,\BBU_n)
\mapsto X_n$ is continuous with respect to the $\rho_\gamma$ topology
(i.e., the rough path topology) for any $\gamma> 1/3$. In particular,
the estimates~\eqref{e:cts_kolm}, combined with the iterated invariance
principle, guarantee that $(U_n,\BBU_n) \to_w (U,\BBU)$ in the $\rho
_\gamma$ topology for some $\gamma> 1/3$. It follows that $X_n \to_w
X$ where $X$ satisfies the rough differential equation
\[
X(t) = X(0) + \int_0^t a\bigl(X(s)\bigr)\,ds +
\int_0^t b\bigl(X(s)\bigr)\,d(U,\BBU) (s).
\]
By definition of rough integrals, and the decomposition \eqref
{e:strato_E}, $X$ satisfies \eqref{e:RDE_SDE}.

Similarly, if the estimates \eqref{e:cts_kolm} hold for all $q <
\infty
$, then we can apply \cite{friz10}, Theorem~12.10, under the relaxed condition
$a\in C^{1+}$, $b\in C^{2+}$.
\end{pf*}

Now let $\phi_t\dvtx\Omega\to\Omega$ be a suspension flow as in
Section~\ref{sec-flow},
with Poincar\'e map $f\dvtx\Lambda\to\Lambda$. As before, we write
$\Omega=\Lambda^r$, $\nu=\mu^r$, where $r\dvtx\Lambda\to\R$ is
a roof
function with
$\bar r=\int r \,d\mu$.
Let $v\dvtx\Omega\to\R^e$ with $\int_\Omega v \,d\nu=0$
and define $\tildev\dvtx\Lambda\to\R^e$,
$\tildev(x)=\int_0^{r(x)}v\circ\phi_t \,dt$.

\begin{cor} \label{cor-approx}
Suppose that $f\dvtx\Lambda\to\Lambda$ is mixing and that $r\in
L^\infty
(\Lambda)$ is bounded away from zero.
Suppose that $a\in C^{1+}$ and $b\in C^3$.
Let $v\in L^\infty(\Omega,\R^e)$ with \mbox{$\int_\Omega v \,d\nu=0$}.
If $\tildev$ admits
a martingale-coboundary decomposition \eqref{eq-decomp_inv} with
$p>\frac{9}2$, then
the conclusion of Theorem~\ref{thmm-approx} is valid.
\end{cor}

\begin{pf}
Recall that $X_n$ satisfies \eqref{e:cts_RDE} with $U_n = W_n$. By
Corollary~\ref{cor-WWflow}, $(W_n,\BBW_n) \to_w (W,\BBW)$ where $W$ is
Brownian motion and
$\BBW^{\beta\gamma}(t) =\break  \int_0^t W^\beta \,dW^\gamma+ D^{\beta
\gamma} t$.
Moreover, the estimates \eqref{e:cts_kolm} hold by Corollary~\ref
{cor-moment_flow}. The result follows directly from Theorem~\ref{thmm:cts_sm}.
\end{pf}

\begin{pf*}{Proof of Theorem~\ref{thmm-approx}}
Again we use the fact that every hyperbolic basic set for an Axiom A
flow can
be written as a suspension over a mixing hyperbolic basic set $f\dvtx
\Lambda
\to\Lambda$ with a H\"older roof function $r$.
Any H\"older mean zero observable $v\dvtx\Lambda\to\R^e$ admits
an $L^\infty$ martingale-coboundary decomposition.
Hence, Theorem~\ref{thmm-approx} follows from Corollary~\ref{cor-approx}
together with the last statement of Theorem~\ref{thmm:cts_sm} (to allow
for the weakened regularity assumption on $b$).
\end{pf*}

\subsection{Rough path theory in discrete time}
In this section, we introduce
tools~\cite{Kelly14} that are the discrete time analogue of those
introduced in the continuous rough path section.
Let $U_n \dvtx[0,T]\to\R^e$ be a step function defined by
\[
U_n(t) = \sum_{j=0}^{[nt]-1}
\Delta U_{n,j}.
\]
We also define the discrete iterated integral $\BBU_n \dvtx[0,T] \to
\R
^{e\times e}$ by
%
\begin{equation}
\label{e:bbu_maps} \BBU_n(t) = \int_0^t
U_n(r) \,dU_n(r) = \sum_{0\leq i < j < [n^{-2}t]}
\Delta U_{n,i} \Delta U_{n,j}.
\end{equation}
Note that, as usual, we use the left-Riemann sum convention. We define
the increments
\[
U_n(s,t) = \sum_{j=[ns]}^{[nt]-1}
\Delta U_{n,j} \quad\mbox{and}\quad \BBU_n(s,t) = \sum
_{[ns] \leq i < j < [n^{-2}t]} \Delta U_{n,i} \Delta U_{n,j}.
\]
Suppose that $a\dvtx\R^d\to\R^d$ is $C^{1+}$ and $b\dvtx\R^d\to
\R^{d\times e}$
is $C^{3}$,
and let $X_{n,j}$ be defined by the recursion
%
\begin{equation}
\label{e:disc_RDE} X_{n,j+1} = X_{n,j} + n^{-1}a(X_{n,j})+
b(X_{n,j})\Delta U_{n,j}
\end{equation}
with initial condition $X_{n,0} = \xi\in\R^d$. We then define the path
$X_n \dvtx[0,T] \to\R^d$ by the rescaling $X_n(t) = X_{n,[nt]}$. The
following theorem is the discrete time analogue of Theorem~\ref
{thmm:cts_sm} and is proved in \cite{Kelly14}.

\begin{thmm}\label{thmm:disc_sm}
Suppose that $(U_n,\BBU_n) \to_w (U,\BBU)$ in $D([0,\infty),\R
^e\times\R
^{e\times e})$, where $U$ is Brownian motion and where $\BBU$ can be written
\[
\BBU(t) = \int_0^t U(s) \,dU(s) + E t,
\]
for some constant matrix $E \in\R^{e\times e}$. Suppose moreover that
there exist $C>0$ and $q>3$ such that
%
\begin{equation}\qquad
\label{e:disc_kolm} \bigl\|U_n(j/n,k/n) \bigr\|_{2q} \leq C\biggl\llvert
\frac{j-k}{n}\biggr\rrvert ^{1/2} \quad\mbox{and} \quad\bigl\|\BBU_n(j/n,k/n)
\bigr\|_{q} \leq C \biggl\llvert \frac
{j-k}{n}\biggr\rrvert ,
\end{equation}
hold for all $n\ge1$ and $j,k = 0, \ldots, n$. Then $X_n \to_w X$ in
$D([0,\infty),\R^d)$, where
\[
dX = \biggl( a(X) + \sum_{\alpha,\beta,\gamma}E^{\beta\gamma
}
\partial ^\alpha b^{\beta}(X) b^{ \alpha\gamma}(X) \biggr) \,dt +
b(X) \,dW.
\]
If \eqref{e:disc_kolm} holds for all $q<\infty$, then the $C^3$
condition on $b$ can be relaxed to $C^{2+}$.
\end{thmm}

\begin{cor} \label{cor-approx2}
Suppose that $f\dvtx\Lambda\to\Lambda$ is mixing and
that $a\in C^{1+}$, $b\in C^3$.
Let $v\in L^\infty(\Omega,\R^e)$ with \mbox{$\int_\Omega v \,d\nu=0$}.
If $v$ admits
a martingale-coboundary decomposition \eqref{eq-decomp_inv} with
$p>\frac{9}2$,
then the conclusion of Theorem~\ref{thmm-approx2} is valid.
\end{cor}

\begin{pf}
We have that $X_n$ is defined by the recursion \eqref{e:disc_RDE} with
$\Delta U_{n,j} = n^{-1/2} v\circ f^j$. In particular, $U_n = W_n$ and
$\BBU_n = \BBW_n$, as defined in Section~\ref{sec-discrete}. By Theorem~\ref{thmm-WW2},
$(W_n,\BBW_n) \to_w (W,\BBW)$ where $W$ is Brownian motion and
$\BBW^{\beta\gamma}(t) = \int_0^t W^\beta \,dW^\gamma+ E^{\beta
\gamma} t$.
Moreover, the estimates \eqref{e:disc_kolm} follow immediately from
Corollary~\ref{cor-moment}. Hence, the result follows from Theorem~\ref
{thmm:disc_sm}.
\end{pf}

\begin{pf*}{Proof of Theorem~\ref{thmm-approx2}}
Again,
any H\"older mean zero observable $v\dvtx\Lambda\to\R^e$ admits
an $L^\infty$ martingale-coboundary decomposition.
Hence, Theorem~\ref{thmm-approx2} follows from Corollary~\ref{cor-approx2}
together with the last statement of Theorem~\ref{thmm:disc_sm}.
\end{pf*}

\section{Generalizations}
\label{sec-gen}

Our main results, Theorems \ref{thmm-WW}, \ref{thmm-approx} for
continuous time,
Theorems \ref{thmm-WW2}, \ref{thmm-approx2} for discrete time,
are formulated for the well known, but restrictive, class of uniformly
hyperbolic (Axiom A) diffeomorphisms and flows.
In this section, we extend these results to a much larger class of systems
that are nonuniformly hyperbolic in the sense of Young \cite{Young98,Young99}.
Also, as promised, we show how to relax the mixing assumption in
Theorems \ref{thmm-WW2}
and \ref{thmm-approx2}.

In Section~\ref{sec-Y1}, we consider the case of noninvertible maps
modeled by Young towers.
Then in Sections \ref{sec-Y2} and \ref{sec-Yflow}, we consider the
corresponding situations for invertible maps and continuous time systems.

\subsection{Noninvertible maps modeled by Young towers}
\label{sec-Y1}

In the noninvertible setting,
a \emph{Young tower} $f\dvtx\Delta\to\Delta$ is defined as follows.
First we recall the notion of a Gibbs--Markov map $F\dvtx Y\to Y$.

Let $(Y,\mu_Y)$ be a probability space with a countable measurable partition
$\alpha$, and let
$F\dvtx Y\to Y$ be a Markov map.
Given $x,y\in Y$, define the separation time $s(x,y)$ to be the least integer
$n\ge0$ such that $F^nx,F^ny$ lie in distinct partition elements of
$\alpha$.
It is assumed that the partition separates orbits.
Given $\theta\in(0,1)$ we define the metric $d_\theta(x,y)=\theta^{s(x,y)}$.

If $v\dvtx Y\to\R$ is measurable, we define $|v|_\theta=\sup_{x\neq
y}|v(x)-v(y)|/d_\theta(x,y)$ and $\|v\|_\theta=\|v\|_\infty
+|v|_\theta$.
The space $F_\theta(Y)$ of observables $v$ with $\|v\|_\theta<\infty$
forms a Banach space with norm $\|\cdot \|_\theta$.

Let $g$ denote the inverse of the Jacobian of $F$ for the measure $\mu
_Y$. We
require the \emph{good distortion} property that $|\log g|_\theta
<\infty$.
The map $F$ is said to be \emph{Gibbs--Markov}
if it has good distortion and \emph{big images}: $\inf_{a\in\alpha
}\mu_Y(Fa)>0$.
A special case of big images is the \emph{full branch} condition $Fa=Y$
for all $a\in\alpha$.
Gibbs--Markov maps with full branches are automatically mixing.

If $F\dvtx Y\to Y$ is a mixing Gibbs--Markov map, then
observables in $F_\theta(Y)$ have exponential decay of correlations
against $L^1$ observables. In particular, Theorems~\ref{thmm-WW2}
and \ref{thmm-approx2} apply in their entirety
to mean zero observables $v\dvtx Y\to\R^e$ with components in
$F_\theta(Y)$
for mixing Gibbs--Markov
maps.

Given a full branch Gibbs--Markov map $F\dvtx Y\to Y$, we now introduce a
\emph{return time function}
$\varphi\dvtx Y\to\Z^+$ assumed to be constant on partition elements.
We suppose that $\varphi$ is integrable and set $\bar\varphi=\int_Y\varphi \,d\mu_Y$.
Define the Young tower
\[
\Delta=\bigl\{(y,\ell)\in Y\times\Z\dvtx0\le\ell<\varphi(y)\bigr\},
\]
and define the tower map
$f\dvtx\Delta\to\Delta$ by setting
%
\begin{equation}
\label{eq-tower} f(y,\ell)= %
\cases{ (y,\ell+1), & \quad $\ell\le\varphi(y)-2$,
\vspace*{2pt}
\cr
(Fy,0), &\quad $\ell=\varphi(y)-1$.}
\end{equation}
Then $\mu=\mu_Y\times{\rm Lebesgue}/\bar\varphi$ is an ergodic
$f$-invariant probability measure on~$\Delta$.
Note that the system $(\Delta,\mu,f)$ is uniquely determined by
$(Y,\mu_Y,F)$ together with~$\varphi$.

The separation time $s(x,y)$ extends to the tower by setting
$s((x,\ell),(y,\ell'))=0$ for $\ell\neq\ell'$ and
$s((x,\ell),(y,\ell))=s(x,y)$. The metric $d_\theta$ extends
accordingly to $\Delta$ and we define the space $F_\theta(\Delta)$ of
observables $v\dvtx\Delta\to\R$ that lie in $L^\infty(\Delta)$
and are
Lipschitz with respect to this metric.

The tower map $f\dvtx\Delta\to\Delta$ is mixing if and only if
$\gcd\{\varphi(a)\dvtx a\in\alpha\}=1$.
In the mixing case, it follows from Young \cite{Young98,Young99} that
the rate of decay of correlations on the tower $\Delta$ is determined
by the tail function
\[
\mu(\varphi>n)=\mu\bigl(y\in Y\dvtx\varphi(y)>n\bigr).
\]
In \cite{Young98}, it is shown that exponential decay of $\mu(\varphi>n)$
implies exponential decay of correlations for observables in $F_\theta
(\Delta)$,
and \cite{Young99} shows that if
$\mu(\varphi>n)=O(n^{-\beta})$ then correlations for such observables
decay at a rate that is
$O(n^{-(\beta-1)})$.
For systems that are \emph{modeled by a Young tower}, H\"older
observables for the underlying dynamical system lift to observables in
$F_\theta(\Delta)$
(for appropriately chosen $\theta$) and thereby inherit the above
results on decay of correlations. Similarly, if we define $F_\theta
(\Delta,\R^e)$ to consist of observables $v\dvtx\Delta\to\R^e$ with
components in $F_\theta(\Delta)$, then
results on weak convergence for vector-valued H\"older observables are
inherited by the lifted observables in
$F_\theta(\Delta,\R^e)$
and so it suffices to prove everything at the Young tower level.

\begin{thmm} \label{thmm-Young}
Suppose that $f\dvtx\Delta\to\Delta$ is a mixing Young tower with return
time function $\varphi\dvtx Y\to\Z^+$ satisfying
$\mu(\varphi>n)=O(n^{-\beta})$.
Let $v\in F_\theta(\Delta,\R^e)$ with $\int_\Delta v \,d\mu=0$. Then:
\begin{longlist}[(a)]
\item[(a)] \textit{Iterated WIP}:
If $\beta>3$, then the conclusions of Theorem~\ref{thmm-WW2} are valid.
\item[(b)] \textit{Convergence to SDE}:
If $\beta>\frac{11}{2}$, then the conclusions of Theorem~\ref
{thmm-approx2} are valid for all $a\in C^{1+}$, $b\in C^{3}$.
\end{longlist}
In particular, Theorems \ref{thmm-WW2} and \ref{thmm-approx2} are
valid for
systems modeled by Young towers with exponential tails
for all $a\in C^{1+}$, $b\in C^{2+}$.
\end{thmm}

\begin{pf}
In the setting of noninvertible (one-sided) Young towers \cite{Young99},
given $v\in F_\theta(\Delta)$ with mean zero, there is a constant $C$
such that
\[
\biggl|\int_\Delta v w\circ f^n \,d\mu \biggr|\le C\|w
\|_\infty n^{-(\beta-1)}\qquad \mbox{for all $w\in L^\infty$, $n
\ge1$.}
\]
Hence, by Proposition~\ref{prop-decomp}, there is an $L^p$
martingale-coboundary decomposition \eqref{eq-decomp} for any $p<\beta-1$.
The desired results follow from Theorem~\ref{thmm-WWmap}
and Corollary~\ref{cor-approx}, respectively.
\end{pf}

If $\beta>2$, or more generally $\varphi\in L^2$, the WIP is well known.
In fact, $\varphi\in L^2$ suffices also for the iterated WIP and the
mixing assumption on $f$ is unnecessary, as shown in the next result.
These assumptions are optimal, since the ordinary CLT is generally
false when $\varphi\notin L^2$.

\begin{thmm} \label{thmm-Young2}
Suppose that $\Delta$ is a Young tower with return time function
$\varphi\in L^2$.
Let $v\in F_\theta(\Delta,\R^e)$ with $\int_\Delta v \,d\mu=0$. Then
$(W_n,\BBW_n)\to_w (W,\BBW)$ where $W$ is an $e$-dimensional
Brownian motion with covariance matrix
\begin{eqnarray*}
\Sigma^{\beta\gamma} & =&\Cov^{\beta\gamma}W(1)
\\
& =&(\bar\varphi)^{-1}\int_Y
\tildev^\beta\tildev^\gamma \,d\mu _Y+(\bar
\varphi)^{-1}\sum_{n=1}^\infty\int
_Y\bigl(\tildev^\beta\tildev^\gamma
\circ F^n+\tildev^\gamma\tildev^\beta\circ
F^n\bigr) \,d\mu_Y,
\end{eqnarray*}
and $\BBW^{\beta\gamma}(t)=\int_0^t W^\beta \,dW^\gamma+E^{\beta
\gamma
}t$ where
\[
E^{\beta\gamma}=(\bar\varphi)^{-1}\sum_{n=1}^\infty
\int_Y\tildev ^\beta \tildev^\gamma
\circ F^n \,d\mu_Y+\int_\Delta
H^\beta v^\gamma \,d\mu,\qquad H(y,\ell)=\sum
_{j=0}^{\ell-1} v(y,j).
\]

If moreover $\mu(\varphi>n)=O(n^{-\beta})$ for some $\beta>\frac{11}{2}$,
then the conclusion of Theorem~\ref{thmm-approx2} (convergence to SDE) holds
for all $a\in C^{1+}$, $b\in C^{3}$.
\end{thmm}

\begin{pf}
We use the discrete analogue of the inducing method used in the proof
of Theorem~\ref{thmm-WWflow}. Define $\tildev\dvtx Y\to\R^e$ by setting
$\tildev(y)=\sum_{j=0}^{\varphi(y)-1}v(f^jy)$.
Then $\tildev$ lies in $L^2$ and $\int_Y \tildev \,d\mu_Y=0$.
Let $P$ denote the transfer operator for $F\dvtx Y\to Y$.
Although $\tildev\notin F_\theta(Y,\R^e)$
an elementary calculation \cite{MN05}, Lemma~2.2, shows that $P\tildev
\in F_\theta(Y,\R^e)$.
In particular, $P\tildev$ has exponential decay of correlations
against $L^1$
observables. It follows that $\chi=\sum_{j=1}^\infty P^j\tildev$
converges in
$L^\infty$, and hence following the proof of Proposition~\ref
{prop-decomp}, we obtain that $\tildev$ admits an $L^2$
martingale-coboundary decomposition.

Define the cadlag processes $\tildeW_n$, $\tildeBBW_n$ as in \eqref
{eq-WWW} using $\tildev$ instead of $v$.
It follows from Theorem~\ref{thmm-WWmap} that
$(\tildeW_n,\tildeBBW_n)\to_w (\tildeW,\tildeBBW)$ where
$\tildeW$ is an $e$-dimensional Brownian motion
and $\tildeBBW^{\beta\gamma}(t)=\int_0^t\tildeW^\beta \,d\tildeW
^\gamma
+\tilde E^{\beta\gamma}t$ with
\[
\Cov^{\beta\gamma}\tilde W(1)=\int_Y
\tildev^\beta\tildev^\gamma \,d\mu _Y+\sum
_{n=1}^\infty\int_Y\bigl(
\tildev^\beta\tildev^\gamma\circ F^n+
\tildev^\gamma\tildev^\beta\circ F^n\bigr) \,d
\mu_Y,
\]
and
$\tilde E^{\beta\gamma}= \sum_{n=1}^\infty\int_Y \tildev^\beta
\tildev
^\gamma\circ F^n \,d\mu_Y$.

Arguing as in the proof of Theorem~\ref{thmm-WWflow}, and noting
Remark~\ref{rmk-reg},
we obtain that $(W_n,\BBW_n)\to_w(W,\BBW)$ where
\begin{eqnarray*}
 W&=&(\bar\varphi)^{-1/2}\tildeW,\qquad \BBW^{\beta\gamma}(t)=\int
_0^t W^\beta \,dW^\gamma+E^{\beta\gamma}t,
\\
E^{\beta\gamma}&=&(\bar\varphi)^{-1}\tilde E^{\beta\gamma
}+\int
_\Delta H^\beta v^\gamma \,d\mu.
\end{eqnarray*}

Finally, to prove the last statement of the theorem, it suffices by
Corollary~\ref{cor-approx} to show that $v$ admits an $L^p$
martingale-coboundary decomposition with $p>\frac{9}2$. We already saw
that this holds for $\Delta$ mixing, equivalently $d=\gcd\{\varphi
(a)\dvtx a\in\alpha\}=1$. If $d>1$, then $\Delta$ can be written as a
disjoint union of $d$ towers $\Delta_k$ each with a Gibbs--Markov map
that is a copy of $F$ and return time function $1_{\Delta_k}\varphi/d$.
Each of these $d$ towers is mixing under
$f^d$, and the towers are cyclically permuted by $f$. Hence,
\[
\sum_{m=1}^\infty P^m\tildev=
\sum_{k,r=1}^\infty\sum
_{m=0}^\infty P^{md+r} \biggl(1_{\Delta_k}
\tildev-d\int_\Delta1_{\Delta_k}\tildev \,d\mu \biggr).
\]
But $\|P^{md}(1_{\Delta_k}\tildev-d\int_\Delta1_{\Delta_k}\tildev
\,d\mu
)\|_p\ll m^{-\beta}$.
Hence, we can define $\chi=\break \sum_{m=1}^\infty P^m\tildev\in
L^p$ yielding the desired decomposition $\tildev=m+\chi\circ f-\chi$.
\end{pf}

\begin{examp}
A prototypical family of nonuniformly expanding maps are intermittent
maps $f\dvtx[0,1]\to[0,1]$ of Pomeau--Manneville type \cite
{PomeauManneville80,LiveraniSaussolVaienti99}
given by
\[
fx= %
\cases{ x\bigl(1+2^\alpha x^\alpha\bigr), &\quad $x
\in\bigl[0,\frac{1}2\bigr),$ \vspace *{2pt}
\cr
2x-1, & \quad $x\in \bigl[
\frac{1}2,1\bigr].$}
\]
For each $\alpha\in[0,1)$, there is a unique absolutely continuous
invariant probability measure $\mu$. For $\alpha\in(0,1)$, there
is a neutral fixed point at $0$ and the system is modeled by a mixing
Young tower with tails that are $O(n^{-\beta})$ where $\beta=\alpha^{-1}$.

Hence, the results of this paper apply in their entirety for $\alpha
\in
[0,\frac{2}{11})$. Further, it is well known that the WIP holds if and
only if $\alpha\in[0,\frac{1}2)$, and we recover this result, together
with the iterated WIP, for $\alpha\in[0,\frac{1}2)$.
\end{examp}

\subsection{Invertible maps modeled by Young towers}
\label{sec-Y2}

A large class of nonuniformly hyperbolic diffeomorphisms (possibly with
singularities) can be modeled by two-sided Young towers with
exponential and polynomial tails. For such towers, Theorems \ref
{thmm-Young} and \ref{thmm-Young2} go through essentially without change.
The definitions are much more technical, but we sketch some of the
details here.

Let $(M,d)$ be a Riemannian manifold.
Young \cite{Young98} introduced a class of nonuniformly hyperbolic maps
$T\dvtx M\to M$ with the property that there is an ergodic
$T$-invariant SRB
measure for which exponential decay of correlations holds for H\"older
observables. We refer to \cite{Young98} for the precise definitions,
and restrict here to providing the notions and notation required for
understanding the results presented here.
In particular, there is a ``uniformly hyperbolic'' subset $Y\subset M$
with partition $\{Y_j\}$ and return time function
$\varphi\dvtx Y\to\Z^+$ (denoted $R$ in \cite{Young98}) constant on
partition elements. For each $j$, it is assumed that
$T^{\varphi(j)}(Y_j)\subset Y$. We define the \emph{induced map}
$F=T^\varphi\dvtx Y\to Y$.

Define the \emph{(two-sided) Young tower} $\Delta=\{(y,\ell)\in
Y\times\Z
\dvtx0\le\ell<\varphi(y)\}$ and define the
tower map $f\dvtx\Delta\to\Delta$ using the formula \eqref{eq-tower}.

It is assumed moreover that there is an $F$-invariant foliation of
$Y$ by ``stable disks,'' and that this foliation
extends up the tower $\Delta$. We obtain the quotient tower map
$\bar f\dvtx\bar\Delta\to\bar\Delta$ and quotient induced map
$\bar F=\bar f^\varphi\dvtx\bar Y\to\bar Y$.
The hypotheses in \cite{Young98} guarantee that:

\begin{prop} \label{prop-Y98}
There exists an ergodic $T$-invariant probability measure $\nu$ on $M$,
and ergodic
invariant probability measures $\mu_\Delta$, $\mu_{\bar\Delta}$,
$\mu_Y$, $\mu_{\bar Y}$ defined on $\Delta$, $\bar\Delta$, $Y$,
$\bar Y$,
respectively, such that:
\begin{longlist}[(a)]
\item[(a)] The projection $\pi\dvtx\Delta\to M$ given by $\pi
(y,\ell
)=T^\ell y$,
and the projections
$\bar\pi\dvtx\Delta\to\bar\Delta$ and $\bar\pi\dvtx Y\to\bar
Y$ given by quotienting,
are measure preserving.
\item[(b)] The return time function $\varphi\dvtx Y\to\Z^+$ is integrable
with respect to $\mu_Y$ (and hence also with respect to $\mu_{\bar Y}$
when regarded as a function on $\bar Y$).
\item[(c)] $\mu_\Delta=\mu_Y\times{\rm counting}/\int_Y\varphi
\,d\mu$ and
$\mu_{\bar\Delta}=\mu_{\bar Y}\times{\rm counting}/\int_Y\varphi
\,d\mu$.
\item[(d)] The system $(\bar Y,\bar F,\mu_{\bar Y})$ is a full branch
Gibbs--Markov map with partition $\alpha=\{\bar Y_j\}$. Hence, $\bar
f\dvtx\bar\Delta\to\bar\Delta$ is a one-sided Young tower as
in Section~\ref{sec-Y1}.
\item[(e)] $\mu_Y(\varphi>n)=O(e^{-an})$ for some $a>0$.
\item[(f)]
Let $v\dvtx M\to\R$ be H\"older with $\int_M v \,d\nu=0$.
Then $v\circ\pi=\bar v\circ\bar\pi+\chi_1\circ f-\chi_1$ where
$\chi
_1\in L^\infty(\Delta)$
and $\bar v\in F_\theta(\bar\Delta)$ for some $\theta\in(0,1)$.
\end{longlist}
\end{prop}

\begin{pf}
Parts (a)--(e) can be found in \cite{Young98}. For part (f) see, for
example, \cite{M07,MN05}.
\end{pf}

\begin{cor} \label{cor-Y98}
Theorems \ref{thmm-WW2} and \ref{thmm-approx2} are valid for H\"older
mean zero observables of systems modeled by (two-sided) mixing Young
towers with exponential tails.
\end{cor}

\begin{pf}
By Proposition~\ref{prop-Y98}(d) and the proof of Theorem~\ref{thmm-Young},
for any $p$ we can decompose $\bar v\in F_\theta(\Delta)$ as
$\bar v=\bar m+\bar\chi_2\circ\bar f-\bar\chi_2$
where $\bar m,\bar\chi_2\in L^\infty(\bar\Delta)$ and $\bar m$
lies in
the kernel of the transfer operator corresponding to $\bar F\dvtx\bar
\Delta
\to\bar\Delta$. Now let $m=\bar m\circ\bar\pi$ and
$\chi=\chi_1+\bar\chi_2\circ\bar\pi$ where $\chi_1$ is as in
Proposition~\ref{prop-Y98}(f).
We have shown that $v\circ\pi$ admits an $L^\infty$ martingale-coboundary
decomposition \eqref{eq-decomp_inv}. By Theorem~\ref{thmm-WWmap_inv}, we
obtain the required results for $v\circ\pi$, and hence for $v$.
\end{pf}

By \cite{BenedicksYoung00}, this includes the important example of H\'
enon-like attractors.
Again the results hold with the appropriate modifications (in the
formulas for $\Sigma$ and $E$)
for nonmixing towers with exponential tails.

A similar situation holds for systems modeled by (two-sided) Young
towers with polynomial tails where Proposition~\ref{prop-Y98}(a)--(d)
are unchanged and part (e) is replaced by the
condition that $\mu_Y(\varphi>n)=O(n^{-\beta})$.
In general, part (f) needs modifying.
The simplest case is where
there is sufficiently fast uniform contraction along stable manifolds
(exponential
as assumed in \cite{AlvesPinheiro08,M07,MN05}, or polynomial as in~\cite
{AlvesAzevedo-sub}).
Then part (f) is unchanged allowing us to reduce to the situations of
Theorem~\ref{thmm-Young} in the mixing case, $\beta>3$, and
Theorem~\ref{thmm-Young2} in the remaining cases.

In the general setting of Young towers with subexponential tails, there
is contraction/expansion only on visits to $Y$ and Proposition~\ref
{prop-Y98}(f) fails.
In this case, an alternative construction \cite{MVarandas} can be used to
reduce from $M$ to $Y$ and then to $\bar Y$.
Define the induced observable $\tilde v$ on $Y$ by setting
$\tilde v(y)=\sum_{\ell=0}^{\varphi(y)-1}v(T^\ell y)$. If $\varphi
\in
L^p$ (which is the case for all $p<\beta$) then it is shown in \cite
{MVarandas}
that $\tilde v=\bar m\circ\bar\pi+\chi\circ F-\chi$
where $\bar m\in L^p(\bar Y)$ lies in the kernel of the transfer
operator for
$\bar F\dvtx\bar Y\to\bar Y$ and $\chi\in L^p(Y)$.
Thus, if $\varphi\in L^2$, we obtain the
iterated WIP for $\tilde v$, and hence for $v$.

\subsection{Semiflows and flows modeled by Young towers}
\label{sec-Yflow}

Finally, we note that the results for noninvertible and invertible maps
modeled by a Young tower pass over to suspension semiflows and flows defined
over such maps.
Using the methods in Sections \ref{sec-flow} and \ref{sec-moment}, we
reduce from observables defined on the flow to observables defined on
the Young tower, where we can apply the results from Sections \ref
{sec-Y1} and \ref{sec-Y2}.
We refer to \cite{MN05} for a description of numerous examples of
flows that
can be reduced to maps in this way.

We mention here the classical Lorenz attractor for which Theorems \ref{thmm-WW}
and \ref{thmm-approx} follow as a consequence of such a construction.
There are numerous methods to proceed with the Lorenz attractor, but
probably the simplest is as follows. The Poincar\'e map is a Young
tower with exponential tails, but the roof function for the flow has a
logarithmic singularity, and hence is unbounded. An idea in \cite
{BMsub} is to remodel the flow as a suspension with bounded roof
function over a mixing Young tower $\Delta$ with slightly worse, namely
stretched exponential, tails. In particular, the return time function
for $\Delta$ still lies
in $L^p$ for all $p$. H\"older observables for the flow can now be
shown to induce to observables
in $F_\theta(\Delta)$, thereby reducing to the situation of Section~\ref{sec-Y2}.
Moreover, the flow for the Lorenz attractor has exponential contraction
along stable manifolds, and this is inherited by each of the Young
tower models described above. Hence, we can reduce to the situation in
Theorem~\ref{thmm-Young} with $\beta$ arbitrarily large.

\section*{Acknowledgements}
We are grateful to Martin Hairer and Andrew Stuart for helpful comments.

%

\printaddresses
\end{document}